\documentclass[10pt]{article}

\title{Partial mean field limits in heterogeneous networks}

\author{
Carsten Chong\thanks{Center for Mathematical Sciences, Technical University of Munich, Boltzmannstra\ss e 3, 85748 Garching, Germany, e-mail: carsten.chong@tum.de, cklu@ma.tum.de, URL: www.statistics.ma.tum.de}
~and
Claudia Kl\"uppelberg$^\ast$
}

\usepackage{latexsym}
\usepackage{amssymb}
\usepackage{amsmath}
\usepackage{amsthm}
\usepackage{mathrsfs}
\usepackage[numbers]{natbib}
\usepackage{url}
\usepackage{color}
\usepackage{dsfont}
\usepackage{longtable}

\usepackage[latin1]{inputenc}
\usepackage{enumitem}
\setenumerate{leftmargin=*}
\usepackage{amsopn}
\usepackage{comment}

\newcommand{\bfi}{\begin{fig}}
\newcommand{\efi}{\end{fig}}
\newcommand{\btab}{\begin{tab}}
\newcommand{\etab}{\end{tab}}
\newcommand{\barr}{\begin{array}}
\newcommand{\earr}{\end{array}}
\newcommand{\beq}{\begin{equation}}
\newcommand{\eeq}{\end{equation}}
\newcommand{\bdis}{\begin{displaymath}}
\newcommand{\edis}{\end{displaymath}\noindent}

\newcommand{\bbn}{\mathbb{N}}

\newcommand{\bbr}{\mathbb{R}}

\newcommand{\bbe}{\mathbb{E}}

\newcommand{\bbp}{\mathbb{P}}

\newcommand{\bbf}{\mathbb{F}}

\newcommand{\bone}{\mathds 1}

\newcommand{\halmos}{\quad\hfill $\Box$}

\newcommand{\calc}{{\cal C}}

\newcommand{\calf}{{\cal F}}

\newcommand{\calp}{{\cal P}}

\newcommand{\Omicron}{{\cal O}}

\newcommand{\calg}{{\cal G}}

\newcommand{\cald}{{\cal D}}

\newcommand{\al}{{\alpha}}
\newcommand{\la}{{\lambda}}
\newcommand{\La}{{\Lambda}}
\newcommand{\eps}{{\epsilon}}

\newcommand{\ga}{{\gamma}}
\newcommand{\Ga}{{\Gamma}}

\newcommand{\si}{{\sigma}}
\newcommand{\Si}{{\Sigma}}

\newcommand{\Om}{{\Omega}}

\newcommand{\var}{{\mathrm{Var}}}
\newcommand{\cov}{{\mathrm{Cov}}}

\newcommand{\dd}{\mathrm{d}}
\newcommand{\ee}{\mathrm{e}}

\newcommand{\PP}{\mathrm{P}}

\newcommand{\CC}{\mathrm{C}}
\newcommand{\inn}{\mathrm{in}}
\newcommand{\out}{\mathrm{out}}

\newcommand{\diag}{\mathrm{diag}}

\makeatletter
\newcommand{\opnorm}{\@ifstar\@opnorms\@opnorm}
\newcommand{\@opnorms}[1]{%
  \left|\mkern-1.5mu\left|\mkern-1.5mu\left|
   #1
  \right|\mkern-1.5mu\right|\mkern-1.5mu\right|
}
\newcommand{\@opnorm}[2][]{%
  \mathopen{#1|\mkern-1.5mu#1|\mkern-1.5mu#1|}
  #2
  \mathclose{#1|\mkern-1.5mu#1|\mkern-1.5mu#1|}
}
\makeatother

\newtheoremstyle{neu}
    {11pt}      
    {11pt}      
    {}                  
    {}          
    {\bfseries} 
    {}          
    {1em}  
    {\textbf{\thmname{#1}\thmnumber{ #2}\thmnote{ (#3)}}}          
    
\newtheoremstyle{proof}
    {11pt}      
    {11pt}      
    {}                  
    {}          
    {\bfseries} 
    {}            
    {1em}          
    {\textbf{\thmname{#1}.}}          

\newtheorem{Theorem}{Theorem}[section]
\newtheorem{Corollary}[Theorem]{Corollary}
\newtheorem{Lemma}[Theorem]{Lemma}
\newtheorem{Proposition}[Theorem]{Proposition}
\theoremstyle{neu}
\newtheorem{Definition}[Theorem]{Definition}
\newtheorem{Example}[Theorem]{Example}
\newtheorem{Remark}[Theorem]{Remark}
\newtheorem{Assumption}{Assumption}

\theoremstyle{proof}
\newtheorem{Proof}{Proof}

\newcommand{\bthm}{\begin{Theorem}}
\newcommand{\ethm}{\end{Theorem}}

\newcommand{\bcor}{\begin{Corollary}}
\newcommand{\ecor}{\end{Corollary}}

\newcommand{\blem}{\begin{Lemma}}
\newcommand{\elem}{\end{Lemma}}

\newcommand{\bprop}{\begin{Proposition}}
\newcommand{\eprop}{\end{Proposition}}

\newcommand{\bdf}{\begin{Definition}}
\newcommand{\edf}{\end{Definition}}

\newcommand{\bex}{\begin{Example}}
\newcommand{\eex}{\end{Example}}

\newcommand{\brem}{\begin{Remark}}
\newcommand{\erem}{\end{Remark}}

\newcommand{\bass}{\begin{Assumption}}
\newcommand{\eass}{\end{Assumption}}

\newcommand{\bpr}{\begin{Proof}}
\newcommand{\epr}{\end{Proof}}

\newcommand{\benu}{\begin{enumerate}}
\newcommand{\eenu}{\end{enumerate}}

\newcommand{\bit}{\begin{itemize}}
\newcommand{\eit}{\end{itemize}}

\newcommand{\bff}{\textbf}

\numberwithin{equation}{section}

\allowdisplaybreaks[2]

\begin{document}

\date{}

\maketitle

\begin{abstract}
We investigate systems of interacting stochastic differential equations with two kinds of heterogeneity: one originating from different weights of the linkages, and one concerning their asymptotic relevance when the system becomes large. To capture these effects we define a partial mean field system, and prove a law of large numbers with explicit bounds on the mean squared error. Furthermore, a large deviation result is established under reasonable assumptions. The theory will be illustrated by several examples: on the one hand, we recover the classical results of chaos propagation for homogeneous systems, and on the other hand, we demonstrate the validity of our assumptions for quite general heterogeneous networks including those arising from preferential attachment random graph models. 
\end{abstract}

\vfill

\noindent
\begin{tabbing}
\> secondary: \,\,\,05C80
\end{tabbing}

\vspace{1cm}


\vspace{0.5cm}

\newpage

\section{Introduction}\label{Sect1}

The application of mean field theory to large systems of stochastic differential equations (SDEs) was initiated by McKean's seminal work \citep{McKean66a, McKean66b, McKean67}. In the classical case, an $N$-dimensional interacting particle system is governed by SDEs of the form
\begin{align} \dd X^N_i(t) &=\frac{1}{N-1}\sum_{j\neq i} \big(X^N_j(t)-X^N_i(t)\big)\,\dd t + \dd B_i(t),\quad t\in\bbr_+, \nonumber\\
X^N_i(0)&=X_i(0),\quad i=1,\ldots,N,\label{McKean}\end{align}
with independent starting random variables $X_i(0)$ and independent Brownian motions $B_i$.
As the number of particles increases, the pair dependencies in this coupled system decrease with order $1/N$ such that a law of large numbers applies (see Theorem 1.4 of \citep{Sznitman91}).
Defining
\begin{align} \dd \bar X^N_i(t) &=\frac{1}{N-1}\sum_{j\neq i} \big(\bbe[\bar X^N_j(t)]-\bar X^N_i(t)\big)\,\dd t + \dd B_i(t),\quad t\in\bbr_+,\nonumber\\ 
\bar X^N_i(0)&=X_i(0),\quad i=1,\ldots,N, \label{McKean-mf}\end{align}
there exists for every $T\in\bbr_+$ a constant $C(T)\in\bbr_+$ independent of $N$ such that 
\beq\label{McKean-mf-dist} 
\sup_{i=1,\ldots,N} \bbe\left[\sup_{t\in[0,T]} |X^N_i(t)-\bar X^N_i(t)|^2\right]^{1/2}\leq \frac{C(T)}{\sqrt{N}}.
\eeq
In other words, in a large system, the behaviour of a fixed number of particles evolving according to \eqref{McKean} is well described by the so-called \emph{mean field system} \eqref{McKean-mf},
where all stochastic processes are stochastically independent, a phenomenon that is called \emph{propagation of chaos}. 
Thus, mean field theory provides a model simplification by reducing a many-body problem as in \eqref{McKean} to a one-body problem as in \eqref{McKean-mf} with explicit $L^2$-estimates on the occurring error. 
Moreover, it can be shown that the empirical measure of the particles satisfies a large deviation principle as $N\to\infty$, see \citep{Dawson87,Leonard95}. 
There exists a huge literature dealing with this or related topics, and we only mention the review papers \citep{Gaertner88, Sznitman91}, where one can also find further references.

The systems \eqref{McKean} and \eqref{McKean-mf} describe statistically equal or exchangeable particles: any permutation of the indices $i\in\{1,\ldots,N\}$ leads to a system with the same distribution (cf. \citep{Vaillancourt88}). In particle physics such an assumption is certainly reasonable and underlies many other similar models of mean field type, see for example the two treatises \citep{Talagrand11, Talagrand11b} for numerous examples.

However, when mean field models are considered in applications other than statistical mechanics, the homogeneity assumption may not be appropriate in all situations.   
For instance, in \citep{Bouchaud00, Ichinomiya12} the processes \eqref{McKean} are used to model the wealth of trading agents in an economy, who are typically far from being equal in their trading behaviour (there are ``market makers'' and others).
Similarly, the stochastic Cucker-Smale model that is considered in \citep{Ahn10, Bolley11} describes the ``flocking'' phenomenon of individuals. 
Also here it only seems natural that one or several ``leaders'' may have a distinguished role, setting them apart from the remaining system. 
Moreover, in systemic risk modelling the particles represent financial institutions that interact with each other through mutual exposures, see \citep{Battiston12, Fouque13, Kley14} for some approaches in this direction.
The different players in the banking sector vary considerably in size and importance, which is obvious from the fact that some banks were considered too big to fail during the financial crisis of 2007-08.
Further fields of applications where mean field theory is used for interacting particle systems include genetic algorithms \citep{DelMoral01}, neuron modelling (see \citep{Fournier15} and references therein) and epidemics modelling \citep{Leonard90}.

Partly triggered by the examples in the previous paragraph, this research aims to investigate deviations from homogeneous systems to heterogeneous systems.
First, we allow for different interaction rates between pairs (instead of $1/(N-1)$ throughout), and second, we permit the subsistence of a core--periphery structure in the mean field limit, that is, some particles may have a non-vanishing influence even when the system becomes large. Another restriction we will relax in our analysis concerns the driving noises of the interacting SDEs: instead of independence we explicitly allow for different degrees of dependence in the noise terms, even asymptotically. Until now there exists only a very small amount of literature that generalizes \eqref{McKean} in these directions: in \citep{Budhiraja15, Kley14, Nagasawa87, Nagasawa87b}  the particles are divided into finitely many groups within which they are homogeneous (and the number of members in both groups must tend to infinity for the law of large numbers), and \citep{Buckdahn14, Carmona14}, where one major agent exists and propagation of chaos for the minor agents is considered conditioned on the major one. 
  Other papers that consider general heterogeneous systems include \citep{Degond14}, where the propagation of chaos result is \emph{assumed}, and \citep{Finnoff93, Finnoff94, Giesecke15}, where a law of large numbers for the empirical measure is proved under various conditions. 
Regarding the last-mentioned papers, two aspects are worth commenting on. First, assuming that finitely many core particles do exist in the system, their contribution to the empirical distribution becomes less and less as $N\to\infty$ although their impact may very well stay high. Thus, in this case the empirical distribution may fail to describe the behaviour of the system as a whole. Second, whereas for homogeneous systems the convergence of the empirical measure is equivalent to the existence of a mean field limit in the sense of \eqref{McKean-mf-dist} (see e.g. Proposition~2.2(i) of \citep{Sznitman91}), this is no longer true for heterogeneous systems. For core particles the left-hand side of \eqref{McKean-mf-dist} need not converge to $0$ even if the empirical distribution converges, say, to a deterministic limit. For example, in the case of \citep{Buckdahn14, Carmona14} with one core particle, an unconditional propagation of chaos result does not hold for this particle without further assumptions (even if it does for the periphery particles). 

Due to the two aforementioned reasons, we will \emph{not} work with the empirical distribution in this paper but state and prove mean field limit theorems for the particles on the process level. In Section~\ref{Sect2} we start by introducing the precise interacting particle model we want to investigate. Then we define a corresponding \emph{partial mean field model}, for which we prove a law of large numbers type result (Theorem~\ref{LLN}) with explicit convergence rates in Section~\ref{Sect3}. It generalizes \eqref{McKean-mf-dist} by taking into account the different kinds of heterogeneity due to varying pair interaction rates, a distinction between important/core and less important/periphery pair relationships, and interdependencies between the driving noise terms. 

The main difficulty here is to identify the correct rates that govern the distance between the original system and the mean field approximation. As we will see, a total of twelve rates is required, each expressing a connectivity property of the underlying interaction and correlation networks. This is inevitable in contrast to \cite{Buckdahn14, Carmona14} where the stochastic dependencies among the particles are annihilated simply by conditioning. In order to elucidate the meaning of each rate, we discuss three exemplary situations in detail. In Section~\ref{Sect31}, in particular in Example~\ref{classex}, we show that in the quasi-homogeneous case all twelve rates typically boil down to a single rate like in \eqref{McKean-mf-dist}. In Section~\ref{Sect32}, we explain why the prerequisites for Theorem~\ref{LLN} in the heterogeneous case are essentially sparsity assumptions on the particle network, which are satisfied for instance if this network is generated from a preferential attachment mechanism, see Section~\ref{Sect33}. In order to show the last statement, we have to derive the asymptotics of the maximal in- and out-degrees of directed preferential attachment graphs, see Lemma~\ref{prefatt}. This result may be of independent interest and generalizes that of \citep{Mori05} for undirected graphs. 

The second main result of our paper is a large deviation principle for the difference $X^N-\bar X^N$, which is presented in Section~\ref{Sect4} as Theorem~\ref{LD}. In contrast to homogeneous systems, where such a principle is proved for the empirical measure (see \citep{Dawson87, Leonard95}), we work on the process level again and therefore need to require the existence of all exponential moments. Furthermore, due to heterogeneity, we do not obtain an explicit formula for the large deviation rate function, but a variational representation as Fenchel-Legendre transform. The final Section~\ref{Sect5} contains the proofs.

\section{The model}\label{Sect2}

Before we introduce the model we analyze in this paper, we list a number of notations that will be employed throughout the paper.
\begin{longtable}{p{3 cm} p{12 cm}}
$\bbr_+$ & the set $[0,\infty)$ of \emph{positive} real numbers;\\
$[z]$ & the largest integer smaller or equal to $z\in\bbr$;\\
$\bbn$ & the natural numbers $\{1, 2, \ldots\}$;\\
$A, x$ & the typical notation for a matrix $A=(A_{ij}\colon i,j\in\bbn)\in\bbr^{\bbn\times\bbn}$ and a vector $x=(x_i\colon i\in\bbn)^\prime\in\bbr^\bbn$, with all binary relations such as $\leq$, or operations relying on them such as the absolute value $|\cdot|$ or taking the supremum being understood componentwise when applied to matrices and vectors;\\
$(\cdot)^\prime$ & the transposition operator;\\
$AB, Ax, \ee^A$ & matrix--matrix and matrix--vector multiplication and the matrix exponential, all defined in analogy to the finite-dimensional case, provided that the involved series converge;\\
$x.y$ & the entrywise product $x.y=(x_iy_i\colon i\in\bbn)^\prime$ for $x,y\in\bbr^\bbn$;\\
$|A|_\infty, |x|_\infty$ & $|A|_\infty:= \sup_{i\in\bbn} \sum_{j\in\bbn} |A_{ij}|$ and $|x|_\infty := \sup_{i\in\bbn} |x_i|$ for $A\in\bbr^{\bbn\times\bbn}$ and  $x\in\bbr^\bbn$;\\
$|A|_\dd$ & $|A|_\dd := \sup_{i\in\bbn} |A_{ii}|$ for matrices $A$;\\
$A^\times$ & the matrix $A$ with all diagonal entries set to $0$;\\
$\mathrm{I}$ & the identity matrix in $\bbr^{\bbn\times\bbn}$ or $\bbr^{d\times d}$ for some $d\in\bbn$;\\
$L^p$ & the space $L^p(\Om,\calf,\bbp)$, $p\in[1,\infty]$, endowed with the topology induced by $\|X\|_{L^p}:=\bbe[|X|^p]^{1/p}$, and to be understood entrywise when applied to matrix- or vector-valued random variables;\\
$\bbe[X], \var[X]$ & componentwise expectation and variance for random variables in $\bbr^{\bbn\times\bbn}$ or $\bbr^\bbn$;\\
$\cov[X,Y],\cov[X]$ & the matrices whose $(ij)$-th entry is $\cov[X_i,Y_j]$ and $\cov[X_i,X_j]$, respectively, when $X$ and $Y$ are random vectors;\\
$x^\ast$ & $x^\ast(t):= \sup_{s\in[0,t]} |x(s)|$ for $t\in\bbr_+$ and functions $x\colon \bbr_+\to\bbr$, again considered entrywise when $x$ takes values in $\bbr^{\bbn\times\bbn}$ or $\bbr^\bbn$;\\
$D^d_T$, $D^\infty_T$ & the space of $\bbr^d$-valued (resp. $\bbr^\bbn$-valued) functions on $[0,T]$ whose coordinates are all c\`adl\`ag functions;\\
$C^d_T, C^\infty_T$ & elements of $D^d_T$ and $D^\infty_T$ where each coordinate is a continuous function;\\
$AC^d_T, AC^\infty_T$ & elements of $D^d_T$ and $D^\infty_T$ where each coordinate is an absolutely continuous function;\\
$\cald^d_T$, $\cald^\infty_T$ & the $\si$-field on $D^d_T$ (resp. $D^\infty_T$) generated by the evaluation maps $\pi_t(x)=x(t)$, $x\in D^d_T$ (resp. $D^\infty_T$), for $t\in[0,T]$;\\
$U,J_1$ & the uniform topology and the Skorokhod topology on $D^d_T$ and $D^\infty_T$ (in the latter case they are defined via the product of the $d$-dimensional topologies);\\
$M^d_T$ & the space of all $(\theta_1,\ldots,\theta_d)$ where each $\theta_i$ is a signed Borel measure on $[0,T]$ of finite total variation $|\theta_i|([0,T])$
\end{longtable}

Given a stochastic basis $(\Om,\calf,\bbf=(\calf(t))_{t\in\bbr_+},\bbp)$ satisfying the usual hypotheses of completeness and right-continuity, we investigate a network described by the following \emph{interacting particle system} (IPS):
\begin{align} \dd X_i(t)&=\sum_{j=1}^\infty a_{ij}(t)X_j(t)\,\dd t + \sum_{j=1}^\infty \si_{ij}(t)X_j(t-)\,\dd L_i(t) + \sum_{j=1}^\infty f_{ij}(t)\,\dd B_j(t)\nonumber\\
&\quad + \sum_{j=1}^\infty \rho_{ij}(t)\,\dd M_j(t),\quad t\in\bbr_+,\quad i\in\bbn, \label{SDE-i} \end{align}
subjected to some $\calf(0)$-measurable $\bbr^\bbn$-valued initial condition $X(0)$. We will also use the more compact form
\beq\label{SDE} \dd X(t) = a(t)X(t)\,\dd t + \si(t)X(t-).\dd L(t) + f(t)\,\dd B(t) + \rho(t)\,\dd M(t),\quad t\in\bbr_+, \eeq
for \eqref{SDE-i}. The ingredients satisfy the following conditions:
\bit
	\item The two measurable functions $t\mapsto a(t)$ and $t\mapsto \si(t)$ are decomposed into $a=a^\CC+a^\PP$ and $\si=\si^\CC + \si^\PP$ such that for all $T\in\bbr_+$ and $i,j\in\bbn$
\beq\label{compatible} A^\diamond_{ij}(T):=\sup_{t\in[0,T]} |a^\diamond_{ij}(t)|<\infty, \quad \Si^\diamond_{ij}(T):=\sup_{t\in[0,T]}|\si^\diamond_{ij}(t)|<\infty,\quad\diamond\in\{\CC,\PP\}. \eeq
We define $A_{ij}(T) := A^\CC_{ij}(T) + A^\PP_{ij}(T)$ and $\Si_{ij}(T):=\Si^\CC_{ij}(T)+\Si^\PP_{ij}(T)$.
	\item $L$ is an $\bbr^\bbn$-valued $\bbf$-L\'evy process (i.e. an $\bbf$-adapted L\'evy process whose increments are independent of the past $\si$-fields in $\bbf$) with finite second moment and mean $0$.
	\item $M$ is an $\bbr^\bbn$-valued square-integrable martingale on any finite time interval, and $B$ is an $\bbr^\bbn$-valued predictable process such that each coordinate process is of locally finite variation. We assume that $B$ and the predictable quadratic variation process $\langle M,M\rangle$ have progressively measurable Lebesgue densities $b\colon \Om\times\bbr_+\to\bbr^\bbn$ and $c\colon \Om\times\bbr_+\to\bbr^{\bbn\times\bbn}$.
	\item $f$ is the sum of two deterministic measurable functions $f^\CC,f^\PP\colon \bbr_+\to\bbr^{\bbn\times\bbn}$, and $\rho$ the sum of two predictable processes $\rho^\CC,\rho^\PP\colon \Om\times\bbr_+\to\bbr^{\bbn\times\bbn}$.
\eit

Of course, the stochastic integrals behind \eqref{SDE} must make sense: each single integral must be well defined \emph{and} the infinite sums must converge in an appropriate sense. A sufficient condition for the existence of the infinite-dimensional integral is the existence of the one-dimensional ones plus the summability of their $L^2$-norms.

Next, we shall explain the rationale behind the IPS model \eqref{SDE} and the specific choices for the involved processes. By the definition given in \eqref{SDE-i}, the processes $(X_i\colon i\in\bbn)^\prime$ are coupled in two ways in general: first, they interact internally with each other through a drift term (determined by $a$) and a volatility term (determined by $\si$ in conjunction with $L$); and second, they are exposed to the same external forces (given by $B$ and $M$), where $f$ and $\rho$ determine the level of influence these noises have on the particles. In particular, by tuning the parameters $a$, $\si$, $f$ and $\rho$ appropriately, one obtains a large range of possible dependence structures for the model \eqref{SDE}.

The question this paper aims to attack is how and to which degree the complexity of the high-dimensional IPS \eqref{SDE} can be reduced. 
Of course, if all entries of the matrices $a$, $\si$, $f$ and $\rho$ are zero or large, there is no hope in simplifying the model. 
Therefore, our focus lies on particle networks, where only a small number of pairs have strong interaction, while the majority of links in the system are relatively weak. This is implemented in the decomposition of $a$, $\si$, $f$ and $\rho$ into a \emph{core matrix} (superscript $\CC$) and a \emph{periphery matrix} part (superscript $\PP$). It is important to notice that our distinction between core and periphery is not made on the basis of the particles, but on the linkages between them. This allows for greater modelling flexibility since it includes multi-tier networks in our analysis.

In the presence of non-negligible pair interactions it is natural to apply the mean field limit only to the links encoded by the periphery matrices.
Therefore, we propose the following \emph{partial mean field system} (PMFS) as an approximation to the IPS \eqref{SDE}:
\begin{align} \dd\bar X(t) &= \Big(a^\CC(t)\bar X(t) + a^\PP(t)\bbe[\bar X(t)] \Big) \,\dd t + \Big(\si^\CC(t)\bar X(t-) + \si^\PP(t)\bbe[\bar X(t)] \Big).\dd L(t)\nonumber\\
&\quad+ f^\CC(t)b(t)\,\dd t + f^\PP(t)\bbe[b(t)]\,\dd t + \rho^\CC(t)\,\dd M(t),\quad t\in\bbr_+,\nonumber\\
 \bar X(0)&=X(0). \label{meanSDE} \end{align}
Written for each row $i\in\bbn$, this is equivalent to:
\begin{align} \dd\bar X_i(t) &= \sum_{j=1}^\infty \Big(a^\CC_{ij}(t)\bar X_j(t) + a^\PP_{ij}(t)\bbe[\bar X_j(t)] \Big) \,\dd t + \sum_{j=1}^\infty \Big(\si^\CC_{ij}(t)\bar X_j(t-) + \si^\PP_{ij}(t)\bbe[\bar X_j(t)] \Big)\,\dd L_i(t)\nonumber\\
&\quad+ \sum_{j=1}^\infty \Big(f^\CC_{ij}(t)b_j(t) +f^\PP_{ij}(t)\bbe[b_j(t)]\Big)\,\dd t + \sum_{j=1}^\infty \rho^\CC_{ij}(t)\,\dd M_j(t),\quad t\in\bbr_+,\nonumber\\
\bar X_i(0)&=X_i(0). \label{meanSDE-i} \end{align}
 
It is clear that a priori there is no reason for \eqref{meanSDE} to be a good approximation for \eqref{SDE}. Therefore, in the next section, we will give precise $L^2$-estimates in terms of the model coefficients for the difference between the IPS and the PMFS. Moreover, we will determine conditions under which this difference becomes small such that we can indeed speak of a law of large numbers. 

 
\section{Law of large numbers}\label{Sect3}

The first main result of this paper assesses the distance between the original IPS \eqref{SDE} and the PMFS \eqref{meanSDE}. To formulate this we have to introduce some further notation. For $T\in\bbr_+$ we define
\begin{align}
v_a(T)&:= \sup_{i\in\bbn} \sum_{j=1}^\infty A_{ij}(T), &v_{a,\dd}(T)&:= \sup_{i\in\bbn} A^\CC_{ii}(T), &v_\si(T)&:=\sup_{i\in\bbn} \sum_{j=1}^\infty \Si_{ij}(T),\nonumber\\
v_L&:=\sup_{i\in\bbn} \|L_i(1)\|_{L^2}, &v_b(T)&:=\sup_{i\in\bbn}\sup_{t\in[0,T]}\|b_i(t)\|_{L^2}, &v_X&:=\sup_{i\in\bbn} \|X_i(0)\|_{L^2},\nonumber\\
v_f(T)&:=\sup_{i\in\bbn} \sup_{t\in[0,T]} \sum_{j=1}^\infty (|f^\CC_{ij}(t)| + |f^\PP_{ij}(t)|),\hspace{-30cm}\nonumber\\
v_{\rho,M}(T)&:=\sup_{i\in\bbn} \sup_{t\in[0,T]} \left(\sum_{j,k=1}^\infty \left|\bbe[\rho^\CC_{ij}(t)\rho^\CC_{ik}(t)c_{jk}(t)]\right|+\left|\bbe[\rho^\PP_{ij}(t)\rho^\PP_{ik}(t)c_{jk}(t)]\right|\right)^{1/2},\hspace{-30cm}\label{finnum}
\end{align}
and introduce the rates
\begin{align} r_1(T)&:=\Big|A^\PP(T)|\cov[X(0)]|(A^\PP(T))^\prime\Big|^{1/2}_\dd, &r_2(T)&:=\Big|\Si^\PP(T)|\cov[X(0)]|(\Si^\PP(T))^\prime\Big|^{1/2}_\dd,\nonumber\\
r_3(T)&:=\Big|A^\PP(T)|\cov[L(1)]|(A^\PP(T))^\prime\Big|_\dd^{1/2}, &r_4(T)&:=\Big|\Si^\PP(T)|\cov[L(1)]|(\Si^\PP(T))^\prime\Big|_\dd^{1/2},\nonumber\\
r_5(T)&:=\sup_{t\in[0,T]} \Big|f^\PP(t)\cov[b(t)](f^\PP(t))^\prime\Big|_\dd^{1/2}, &r_6(T)&:=\sup_{t\in[0,T]} \Big| \bbe\left[\rho^\PP(t)c(t)(\rho^\PP(t))^\prime\right]\Big|_\dd^{1/2},\nonumber\\
r_7(T)&:=\big|A^\PP(T)A^\CC(T)^\times\big|_\infty, &r_8(T)&:=\big|\Si^\PP(T)A^\CC(T)^\times\big|_\infty,\nonumber\\
r_9(T) &:=\sup_{s,t\in[0,T]}\Big|A^\PP(T)|f^\CC(s)\cov[b(s),b(t)](f^\CC(t))^\prime|(A^\PP(T))^\prime \Big|_\dd^{1/2},\hspace{-30cm}\nonumber\\
r_{10}(T)&:=\sup_{s,t\in[0,T]}\Big|\Si^\PP(T)|f^\CC(s)\cov[b(s),b(t)]f^\CC(t))^\prime|(\Si^\PP(T))^\prime \Big|_\dd^{1/2},\hspace{-30cm}\nonumber\\
r_{11}(T)&:=\sup_{t\in[0,T]}\Big| A^\PP(T)|\bbe[\rho^\CC(t)c(t)(\rho^\CC(t))^\prime]| (A^\PP(T))^\prime\Big|_\dd^{1/2},\hspace{-30cm}\nonumber\\
r_{12}(T)&:=\sup_{t\in[0,T]}\Big| \Si^\PP(T)|\bbe[\rho^\CC(t)c(t)(\rho^\CC(t))^\prime]| (\Si^\PP(T))^\prime\Big|_\dd^{1/2}.\hspace{-30cm} \label{rates}
\end{align}

\bthm\label{LLN} Fix some $T\in\bbr_+$ and grant the general model assumptions as given in Section~\ref{Sect2}. Furthermore, assume that each of the numbers in \eqref{finnum} is finite. Then \eqref{SDE} and \eqref{meanSDE} have a pathwise unique solution $X$ and $\bar X$, respectively, and there exist constants $K(T)$ and $K_\iota(T)$, $\iota=1,\ldots,12$, which depend on the model coefficients only through the numbers in \eqref{finnum}, such that
\begin{align} \sup_{i\in\bbn} \left\|(X_i-\bar X_i)^\ast(T)\right\|_{L^2} &\leq K(T) \sum_{\iota=1}^{12} K_\iota(T)r_\iota(T). \label{LLN1}
\end{align}
\ethm

The proof of Theorem~\ref{LLN} will be given in Section~\ref{Sect5}. Compared to the homogeneous case of \citep{Sznitman91}, we have to take care of several kinds of heterogeneous dependencies in the system: different weights on the edges, the distinction between core and periphery links, and possibly dependent driving noises. This explains why we have twelve rates in contrast to a single one in \eqref{McKean-mf-dist}.

\brem Our calculations furnish the following constants in \eqref{LLN1}:
\beq	K(T):=\sqrt{2}\exp((T^{1/2}v_a(T)+2v_\si(T)v_L)^2 T), \label{KT}\eeq
and
\begin{align*}
K_1(T)&:=E(T)T, &K_2(T)&:=2v_L E(T) T^{1/2},\\
K_3(T)&:=\frac{2}{3}E(T)v_\si(T)V(T)T^{3/2}, &K_4(T)&:=\sqrt{2}v_L E(T)v_\si(T)V(T)T,\\
K_5(T)&:=T, &K_6(T)&:=2T^{1/2}, \\
K_7(T)&:=\frac{1}{2}E(T)V(T)T^2, &K_8(T)&:=\frac{2}{\sqrt{3}} v_L E(T)V(T)\\
K_9(T)&:=\frac{1}{2}E(T)T^2, &K_{10}(T)&:=\frac{2}{\sqrt{3}}v_L E(T) T^{3/2},\\
K_{11}(T)&:=\frac{2}{3}E(T)T^{3/2}, &K_{12}(T)&:=\sqrt{2}v_L E(T) T, 
\end{align*}
where
\[
	E(T):=\ee^{v_{a,\dd}(T)}, \quad 	V(T):=\sqrt{2}\ee^{(v_a(T)T^{1/2} + 2v_L v_\si(T))^2 T} \left(v_X + v_f(T)v_b(T) T + 2v_{\rho,M}(T)T^{1/2}\right). 
\]\halmos
\erem


\brem\label{generalize} There are several possibilities to extend Theorem~\ref{LLN} without substantially new arguments.
\benu
	\item It is straightforward to show that Theorem~\ref{LLN} can be extended to the case where the interaction matrices $a$ and $\si$ are replaced by (still deterministic but possibly history-dependent) linear functionals.
\item Suppose that $L=\Gamma L^0$ with some matrix $\Gamma\in\bbr^{\bbn\times\bbn}$ and some other L\'evy process $L^0$ with finite variance and mean zero. Furthermore, $\Ga=\Ga^\CC + \Ga^\PP$ and accordingly $L^\CC=\Ga^\CC L^0$ and $L^\PP=\Ga^\PP L^0$. What one would like to do when passing to the PMFS \eqref{meanSDE} is to replace $L$ there by $L^\CC$. How does this affect the estimate \eqref{LLN1} in Theorem~\ref{LLN}? A similar analysis as for Theorem~\ref{LLN} reveals that an extra rate 
\[ r_{13}:=\Big|\Ga^\PP\cov[\tilde L(1)](\Ga^\PP)^\prime\Big|_\dd^{1/2} \]
appears with constant $K_{13}:=2v_\si(T)V(T)T^{1/2}$.
\item Two further generalizations are discussed in Remark~\ref{gen2} and  Remark~\ref{gen3} below. \halmos
\eenu
\erem

It is obvious that the usefulness of Theorem~\ref{LLN} depends on the sizes of the rates in \eqref{rates}: only if they are small, the PMFS \eqref{meanSDE} is a good approximation to the IPS \eqref{SDE}. Moreover, there are two different views on Theorem~\ref{LLN}: first, if we assume that the underlying network of the IPS is static, it gives an upper bound on the $L^2$-error when the  IPS is approximated by the PMFS; and second, if the interaction network (i.e. $a$, $\si$, $f$ and $\rho$) is assumed to evolve according to an index $N\in\bbn$, Theorem~\ref{LLN} gives conditions under which the PMFS converges in the $L^2$-sense to the IPS when $N\to\infty$ (this happens precisely when all rates in \eqref{rates} converge to $0$ as $N\to\infty$, and the numbers in \eqref{finnum} are majorized independently of $N$). It is also this second point of view that is the traditional one in mean field analysis and that justifies the title ``Law of large numbers'' for the current section.

In the following subsections we will study three examples of dynamical networks and the corresponding conditions for the law of large numbers to hold for the PMFS.

\subsection{Propagation of chaos}\label{Sect31} We first discuss the phenomenon of chaos propagation, and our results will particularly extend the results of \citep{Fouque13}, Section~17.3, \citep{Kley14}, Corollary~4.1, and \citep{Sznitman91}, Theorem~1.4 by including inhomogeneous weights in the model. The setting is as follows:
\benu
\item The underlying network changes with $N\in\bbn$. In particular, we will index $X$ and $\bar X$, the coefficients $a$, $\si$, $f$ and $\rho$ as well as the rates in \eqref{rates} by $N$.
\item All structural assumptions in Section~\ref{Sect2} hold and the numbers in \eqref{finnum}, some of which now depend on $N$, are uniformly bounded in $N$.
\item The core matrices $a^{N,\CC}(t)$, $\si^{N,\CC}(t)$, $f^{N,\CC}(t)$ and $\rho^{N,\CC}(t)$ are diagonal matrices for all times $t\in\bbr_+$.
\item For each $N\in\bbn$, $(L_i, b_i, M_i, \rho^{N,\CC}_{ii},X^N_i(0)\colon i\in\bbn)$ is a sequence of independent random elements (note that the noises indexed by a fixed $i$ may depend on each other).
\item For each $T\in\bbr_+$ the following rates converge to $0$ as $N\to\infty$:
\begin{align*}
r^N_a(T)&:=\sup_{i\in\bbn} \left(\sum_{j=1}^\infty (A^{N,\PP}_{ij}(T))^2\right)^{1/2}, &r^N_{\rho,M}(T)&:= \sup_{i\in\bbn}\sup_{t\in[0,T]} \left(\sum_{j=1}^\infty \bbe[(\rho^{N,\PP}_{ij}(t))^2c_{jj}(t)]\right)^{1/2},\nonumber\\
r^N_\si(T)&:= \sup_{i\in\bbn} \left(\sum_{j=1}^\infty (\Si^{N,\PP}_{ij}(T))^2\right)^{1/2}, &r^N_f(T)&:=\sup_{i\in\bbn} \sup_{t\in[0,T]} \left(\sum_{j=1}^\infty (f^{N,\PP}_{ij}(t))^2\right)^{1/2}.  
\end{align*}
\eenu
These hypotheses ensure that all pair dependencies between the processes $X^N_i$, $i\in\bbn$, vanish when $N\to\infty$. As a result, in the PMFS, the independence of the particles $i$ at $t=0$ propagates through all times $t>0$: the PMFS decouples in contrast to the original IPS. 

\bex\label{classex} In classical mean field theory as in the references mentioned in the introduction, the $N$-th network consists of exactly $N$ particles. In other words, $a^N_{ij}$,  $\si^N_{ij}$, $f^N_{ij}$, $\rho^N_{ij}$ and $X^N_i(0)$ are all $0$ for $i>N$ or $j>N$. Moreover, all pair interaction is assumed to be of order $1/N$, that is, we have for each $T\in\bbr_+$
\beq A^{N,\PP}_{ij}(T)=\frac{A_{ij}(T)}{N}, \quad\Si^{N,\PP}_{ij}(T)=\frac{\Si_{ij}(T)}{N},\quad i,j\in\bbn, \label{order1N} \eeq
where $A_{ij}(T),\Si_{ij}(T) \in\bbr_+$ are uniformly bounded in $i,j\in\bbn$.
Furthermore, the driving noises are supposed to be independent for different particles and to enter the PMFS completely. This means that (3) and (4) hold and that $f^{N,\PP}=\rho^{N,\PP}=0$.
It is easily shown that under these specifications the rates in (5) above converge to $0$ as $N\to\infty$: $r^N_{\rho,M}(T)$ and $r^N_f(T)$ are simply $0$, and $r^N_a(T)$ and $r^N_\si(T)$ are of order $1/\sqrt{N}$ as $N\to\infty$. \halmos
\eex

We still need to show that under assumptions (1)--(5) above, all rates $r^N_\iota(T)$, $\iota=1,\ldots, 12$, converge to $0$ as $N\to\infty$. Since $A^{N,\CC}(T)$ is diagonal, we have $A^{N,\CC}(T)^\times=0$, and since the driving noises for different particles are independent, all covariances (or covariations) vanish outside the diagonal. Thus, we have
\begin{align*} r^N_1(T) &\leq v_X r^N_a(T), & r^N_2(T) &\leq v_X r^N_\si(T), & r^N_3(T) &\leq v_L r^N_a(T), \\
r^N_4(T) &\leq v_L r^N_\si(T), & r^N_5(T) &\leq v_b(T) r^N_f(T), & r^N_6(T) &=r^N_{\rho,M}(T),\\
r^N_7(T) &= 0, & r^N_8(T) &= 0, & r^N_9(T) &\leq v_b(T) v_f(T) r^N_a(T),\\
r^N_{10}(T)&\leq v_b(T) v_f(T) r^N_\si(T), & r^N_{11}(T)&\leq v_{\rho,M}(T) r^N_a(T), & r^N_{12}(T) &\leq v_{\rho,M}(T) r^N_\si(T),
\end{align*}
which all converge to $0$ as $N\to\infty$ by hypothesis. The following remark continues Remark~\ref{generalize} regarding further generalizations of Theorem~\ref{LLN}.
\brem\label{gen2}
In the setting of this subsection there are actually no core relationships between different particles: every pair interaction rate tends to $0$ with large $N$. If we even assume that there is no dependence at all originating from the noises (i.e. $f^{N,\PP}=\rho^{N,\PP}=0$ above), the propagation of chaos result can easily be extended to nonlinear Lipschitz interaction terms (suitably bounded in $N$) instead of the matrices $a^N$ and $\si^N$. As a matter of fact, the classical method of \citep{Sznitman91}, Theorem~1.4, can be applied with obvious changes. \halmos
\erem

\subsection{Sparse interaction versus sparse correlation}\label{Sect32} The propagation of chaos result in the last subsection was based on two core hypotheses: asymptotically vanishing pair interaction rates and the independence of the particles' driving noises. The motivation for establishing Theorem~\ref{LLN}, however, is to deal with situations where these two conditions are precisely not satisfied, that is, when the coefficients $a$, $\si$, $f$ and $\rho$ of \eqref{SDE} are decomposed into a core and a periphery part in a non-trivial way. In fact, in this subsection we discuss a typical situation where the full generality of Theorem~\ref{LLN} is required. Before that, we recall that we consider networks indexed by $N\in\bbn$, and that we are interested in the cases when the rates in \eqref{rates} vanish when $N$ becomes large. 

\subsubsection*{General assumptions}
The following list of hypotheses describes the setting in this subsection.
\benu
\item The statements (1) and (2) of Section~\ref{Sect31} hold.
\item 
$M$ is an $\bbr^\bbn$-valued $\bbf$-L\'evy process, implying that $c_{ij}(t)=\cov[M_i(1),M_j(1)] t$. 
\item At stage $N$, the system consists of $N_0+N$ particles with some fixed $N_0\in\bbn$, that is, we have $a^N_{ij}=\si^N_{ij}=f^N_{ij}=\rho^N_{ij}=X^N_i(0)=0$ as soon as $i>N_0+N$ or $j>N_0+N$.
\item $\calc:=\{1,\ldots,N_0\}$ contains the \emph{core particles}, $\calp^N:=\{N_0+1,\ldots,N\}$ the \emph{periphery particles}, whose number increases with $N$. 
Correspondingly, $a^{N,\CC}$ and $\si^{N,\CC}$ (resp. $a^{N,\PP}$ and $\si^{N,\PP}$) characterize the influence of the core (resp. periphery) particles in the system. In other words, $j\in\calc$ implies that $a^{N,\PP}_{ij}(t)=\si^{N,\PP}_{ij}(t)=0$ for all $i\in\bbn$ and $t\in\bbr_+$, while $j\in\calp^N$ implies $a^{N,\CC}_{ij}(t)=\si^{N,\CC}_{ij}(t)=0$ for all $i\neq j$ and $t\in\bbr_+$. We assume that the diagonals of $a^N$ and $\si^N$ are completely contained in $a^{N,\CC}$ and $\si^{N,\CC}$, respectively. It follows that the partitions of $a^N$ and $\si^N$ can be illustrated as (omitting all zero rows and columns, and using $\ast$ for all potentially non-zero elements):
\[ \begin{array}{c} a^{N,\CC}/ \\ \si^{N,\CC}\end{array}  = \stackrel{\quad\quad\quad\quad\quad~ N_0\quad\quad\quad\quad\quad\quad\quad N}{\left(\begin{array}{c c c | c c c c} \ast & \cdots & \ast & 0 & \cdots & \cdots & 0 \\ \vdots & \ddots & \vdots & \vdots & \ddots & \ddots & \vdots \\ \ast & \cdots & \ast & 0 & \cdots & \cdots & 0\\ \hline \ast & \cdots & \ast & \ast & 0 & \cdots & 0 \\ \vdots & \ddots & \vdots & 0 & \ddots & \ddots & \vdots \\ \vdots & \ddots & \vdots & \vdots & \ddots & \ddots & 0 \\ \ast & \cdots & \ast & 0 & \cdots & 0 & \ast \end{array}\right)}, \quad\begin{array}{c} a^{N,\PP}/ \\ \si^{N,\PP}\end{array} = \stackrel{\quad\quad\quad\quad\quad~ N_0\quad\quad\quad\quad\quad\quad\quad N}{\left(\begin{array}{c c c | c c c c} 0 & \cdots & 0 & \ast & \cdots & \cdots & \ast \\ \vdots & \ddots & \vdots & \vdots & \ddots & \ddots & \vdots \\ 0 & \cdots & 0 & \ast & \cdots & \cdots & \ast\\ \hline 0 & \cdots & 0 & 0 & \ast & \cdots & \ast \\ \vdots & \ddots & \vdots & \ast & \ddots & \ddots & \vdots \\ \vdots & \ddots & \vdots & \vdots & \ddots & \ddots & \ast \\ 0 & \cdots & 0 & \ast & \cdots & \ast & 0 \end{array}\right)}.\]
\item There is a finite number of \emph{systematic noises}, namely $B_1,\ldots B_{N_{00}}$ and $M_1,\ldots,M_{N_{00}}$ for some fixed $N_{00}\in\bbn$ independent of $N$, that are important to a large part of the system, and there are \emph{idiosyncratic noises} $B_{N_{00}+i}$ and $M_{N_{00}+i}$ that only affect the specific particle $i\in\{1,\ldots,N\}$. Thus, we assume for all $i=1,\ldots,N$ and $t\in\bbr_+$ that $\rho^{N,\PP}_{ij}(t)=f^{N,\PP}_{ij}(t)=0$ for $j\in\{1,\ldots,N_{00}\}\cup\{N_{00}+i\}$ and $\rho^{N,\CC}_{ij}(t)=f^{N,\CC}_{ij}(t)=0$ for the other values of $j$. Hence, $f^N$ and $\rho^N$ are of the form
\[
\begin{array}{c} f^{N,\CC}/ \\ \rho^{N,\CC}\end{array}=\stackrel{\quad\quad\quad\quad\quad\quad N_{00}\quad\quad\quad\quad\quad\quad N_{00}+N}{\left(\begin{array}{c c c | c c c c} \ast & \cdots & \ast & \ast & 0 & \cdots & 0\\ \vdots & \ddots & \vdots & 0 & \ddots & \ddots & 0 \\ \vdots & \ddots & \vdots & \vdots & \ddots & \ddots & 0 \\ \ast & \cdots & \ast & 0 & \cdots & 0 & \ast \end{array}\right)}, \quad\begin{array}{c} f^{N,\PP}/ \\ \rho^{N,\PP}\end{array}=\stackrel{\quad\quad\quad\quad\quad\quad N_{00}\quad\quad\quad\quad\quad\quad N_{00}+N}{\left(\begin{array}{c c c | c c c c} 0 & \cdots & 0 & 0 & \ast & \cdots & \ast \\ \vdots & \ddots & \vdots & \ast & \ddots & \ddots & \vdots \\ \vdots & \ddots & \vdots & \vdots & \ddots & \ddots & \ast \\ 0 & \cdots & 0 & \ast & \cdots & \ast & 0 \end{array}\right)}.
\]
\item We have for all $T\in\bbr_+$
\beq\label{Aspecial} A^{N,\PP}_{ij}(T) = \frac{\phi^N_{ij}(T)}{R^N_A},\quad \Si^{N,\PP}_{ij}(T)=\frac{\psi^N_{ij}(T)}{R^N_\Si},\quad i,j=1,\ldots,N, \eeq
where the rates $R^N_A,R^N_\Si\in\bbr_+$ satisfy
\beq\label{sqrtrate} \frac{R^N_A}{\sqrt{N}}\to\infty, \quad \frac{R^N_\Si}{\sqrt{N}}\to\infty, \quad \text{as }N\to\infty, \eeq
and the numbers $\phi^N_{ij}(T),\psi^N_{ij}(T)\in\bbr_+$ satisfy
\[ \phi(T):=\sup_{i,j,N \in\bbn} \phi^N_{ij}(T) < \infty,\quad \psi(T):=\sup_{i,j,N \in\bbn} \psi^N_{ij}(T) < \infty.\]
Note that we always have $\phi^N_{ii}(T)=\psi^N_{ii}(T)=0$.
\item For different $i,j\in\bbn$, the noises $M_i$ and $M_j$ as well as $B_i$ and $B_j$ are uncorrelated.
\item The rates $r^N_f(T)$ and $r^N_{\rho,M}(T)$ from Section~\ref{Sect31} converge to $0$ as $N\to\infty$ for all $T\in\bbr_+$.
\item For each $N\in\bbn$, the initial values $(X^N_i(0)\colon i\in\calp^N)$ are mutually uncorrelated.
\eenu

Conditions (4) and (5) determine the core--periphery structure of the IPS. 
In practice, 
a fixed distinction between core and periphery particles is often not possible because a large number of particles may be engaged in some strong and some weak linkages at the same time. 
As already pointed out, this does not affect the applicability of Theorem~\ref{LLN}, since the concept of core and periphery refers to the linkages there. The choice of fixed core and periphery \emph{particles} in this subsection is only a special case thereof, intended to simplify the arguments below. Next, regarding (6), one can take $R^N_A,R^N_\Si=N$ for concreteness, which can then be compared with Section~\ref{Sect31}. Furthermore, let us point out that assumption (7) is only for convenience (namely that $f^N$ and $\rho^N$ carry the whole correlation structure of the noises). Indeed, it is always possible (under our second-moment conditions) to replace any stochastic integral $\rho\cdot M$, where $M$ is a L\'evy process, with an arbitrary correlation structure by $\rho^\prime\cdot M^\prime$ where $M^\prime$ consists of mutually uncorrelated L\'evy processes (of course, (8) would change accordingly). Finally, if $X^N(0)$ is independent of the driving noises, (9) can be enforced simply by switching to the conditional distribution given $X^N(0)$. 

Under (1)--(9) it is easy to prove that the rates $r^N_1(T)$, $r^N_2(T)$, $r^N_5(T)$ and $r^N_6(T)$ converge to $0$ when $N\to\infty$. For the latter two, this can be deduced in the same way as in Section~\ref{Sect31} because the driving noises of different particles are uncorrelated. For the other two, we use that the starting random variables of periphery particles are assumed to be uncorrelated. Hence, we have by \eqref{sqrtrate}, as $N\to\infty$, that
	\begin{align*} r^N_1(T) &= \sup_{i\in\bbn}\left(\sum_{j\in\calp^N} (A^{N,\PP}_{ij}(T))^2\var[X^N_j(0)]\right)^{1/2}\leq \phi(T) v_X \frac{\sqrt{N}}{R^N_A}\to0, \\
	r^N_2(T)&= \sup_{i\in\bbn}\left(\sum_{j\in\calp^N} (\Si^{N,\PP}_{ij}(T))^2\var[X^N_j(0)]\right)^{1/2}\leq \psi(T) v_X\frac{\sqrt{N}}{R^N_\Si}\to0.
	\end{align*}

However, the nine conditions above are in general \emph{not} sufficient to imply the smallness of the other rates in \eqref{rates}. We need to add extra hypotheses.

\subsubsection*{Sparseness assumptions}
For each of the remaining rates, we further examine what type of conditions are needed to make them asymptotically small. As we shall see, it is always a mixture of a sparseness condition on the interaction matrices $A^N$ and $\Si^N$ and a sparseness condition on the correlation matrices $f^N$ and $\rho^N$.

\vspace{\baselineskip}
\noindent \underline{$r^N_3(T)$ and $r^N_4(T)$}:\quad We first present a counterexample to show that we have to require further conditions. Consider the simple case where $L_i=L_1$ for all $i\in\bbn$ and that $A^{N,\PP}_{ij}(T)=1/R^N_A$ for all $T\in\bbr_+$ and $i,j\in\{1,\ldots,N_0+N\}$ with $i\neq j$. Then 
	\[ r^N_3(T) =\sup_{i\in\bbn} \left(\sum_{j,k\in\calp^N\setminus\{i\}} \left(\frac{1}{\rho_A^N}\right)^2 \cov[L_j(1),L_k(1)]\right)^{1/2} = v_L\frac{N}{R^N_A}, \]
	which need not to converge to $0$ in general. A similar calculation can be done for $r^N_4(T)$. In order to make the rates $r^N_3(T)$ and $r^N_4(T)$ small, there are basically two options: we require the interaction matrices $A^{N,\PP}$ and $\Si^{N,\PP}$ to be sparse, or we require the correlation matrix of $L$ to be sparse. Any other possibility is a suitable combination of these two.
	\benu
		\item[(10a)] The noises $(L_i\colon i\in\calp^N)$ corresponding to periphery particles only have sparse correlation (which, in particular, includes the case of mutual independence as in Section~\ref{Sect31}). More precisely, we require
		\beq\label{L-uncorr} p^N_L := \#\{(i,j)\in\calp^N\times\calp^N\colon \cov[L_i(1),L_j(1)]\neq0\}| = o\big((R^N_A)^2\wedge(R^N_\Si)^2\big) \eeq
for large $N$. Then
		\[ r^N_3(T)=\sup_{i\in\bbn} \left(\sum_{j,k\in\calp^N} A^{N,\PP}_{ij}(T)A^{N,\PP}_{ik}(T)\cov[L_j(1),L_k(1)]\right)^{1/2} \leq \phi(T) v_L \frac{\sqrt{p^N_L}}{R^N_A}\to0, \]
		and, similarly, $r^N_4(T)\to0$ as $N\to\infty$.
		\item[(10b)] The matrices $A^{N,\PP}(T)$ and $\Si^{N,\PP}(T)$, which describe the influence of periphery particles on the system, are only sparsely occupied, in the sense that every particle in the system is only affected by a small number of periphery particles. 
		In mathematical terms this condition reads as
		\begin{align} p^N_{A,1}(T) &:= \sup_{i\in\bbn} \#\{j\in \calp^N\colon A^{N,\PP}_{ij}(T)\neq0\}=o(R^N_A),\nonumber\\
		p^N_\Si(T) &:= \sup_{i\in\bbn} \#\{j\in \calp^N\colon \Si^{N,\PP}_{ij}(T)\neq0\}=o(R^N_\Si).\label{A-sparse}\end{align}
		In this case, we get
		\[ r^N_3(T)=\sup_{i\in\bbn} \left(\sum_{j,k\in\calp^N} A^{N,\PP}_{ij}(T)A^{N,\PP}_{ik}(T)\cov[L_j(1),L_k(1)]\right)^{1/2} \leq\phi(T) v_L \frac{p^N_{A,1}}{R^N_A} \to0, \]
		and similarly $r^N_4(T)\to0$ as $N\to\infty$.
	\eenu
	
\vspace{\baselineskip}	
\noindent \underline{$r^N_7(T)$ and $r^N_8(T)$}:\quad These two rates express the connectivity between core and periphery particles. In general, they will be not become small with large $N$. For instance, if $A^{N,\CC}_{ij}(T)=1$ for all $j\in\calc$ and $i\neq j$, and $A^{N,\PP}_{ij}=1/R^N_A$ for all $j\in\calp^N$ and $i\neq j$, then
	\[ r^N_7(T) = \sup_{i\in\bbn} \sum_{j\in\calp^N} \sum_{k\in \calc} A^{N,\PP}_{ij}(T) A^{N,\CC}_{jk}(T) = N_0\frac{N}{R^N_A}, \]
does not necessarily converge to $0$. An analogous statement holds for $r^N_8(T)$. For $r^N_7(T), r^N_8(T)\to 0$ we have to require that the lower left block of $A^{N,\CC}$, which describes the influence of core particles on periphery particles, or the matrices $A^{N,\PP}(T)$ and $\Si^{N,\PP}(T)$, which describe the influence of periphery particles on the system, be sparse (or a combination thereof):
	\benu 
		\item[(11a)] The influence of core on periphery particles is sparse. 
		In other words, we suppose for the maximal number of periphery particles a single core particle interacts with through the drift:
		\beq\label{A1-sparse} p^N_{A,2} := \sup_{j\in\calc} \#\{i\in\calp^N\colon A^{N,\CC}_{ij}(T)\neq0\} = o(R^N_A\wedge R^N_\Si). \eeq
		Then,
		\[ r^N_7(T)=\sup_{i\in\bbn} \sum_{j\in\calp^N} \sum_{k\in \calc} A^{N,\PP}_{ij}(T) A^{N,\CC}_{jk}(T) \leq N_0 \phi(T) v_a(T) \frac{p^N_{A,2}}{R^N_A} \to 0 \]
		as well as $r^N_8(T)\to0$ as $N\to\infty$.
		\item[(11b)] $A^{N,\PP}(T)$ and $\Si^{N,\PP}(T)$ are sparse in the sense of \eqref{A-sparse}. Then $r^N_7(T),r^N_8(T)\to0$ follow similarly.
	\eenu

\vspace{\baselineskip}
\noindent \underline{$r^N_9(T)$, $r^N_{10}(T)$, $r^N_{11}(T)$ and $r^N_{12}(T)$}:\quad Similar considerations as before show that these four rates do not converge to $0$ in general. 
Instead, we again need to require some mixture of sparsely correlated driving noises and sparsely occupied matrices $A^{N,\PP}$ and $\Si^{N,\PP}$:
		\benu
			\item[(12a)] 
			We assume that for all $T\in\bbr_+$
			\begin{align}
				\label{pNf} p^N_f(T):=\sup_{j\in\{1,\ldots,N_{00}\}} \#\{i\in \calp^N\colon f^{N,\CC}_{ij}\not\equiv0\text{ on } [0,T]\}=o(R^N_A\wedge R^N_\Si),\\
				\label{pNrho} p^N_\rho(T):=\sup_{j\in\{1,\ldots,N_{00}\}} \#\{i\in \calp^N\colon \rho^{N,\CC}_{ij}\not\equiv0\text{ on } [0,T]\}=o(R^N_A\wedge R^N_\Si).
			\end{align}
				Then, recalling that the components of $b$ and $M$ are mutually uncorrelated,
				\begin{align*} r^N_9(T)&=\sup_{i\in\bbn}\sup_{s,t\in[0,T]} \Bigg(\sum_{j,k\in\calp^N}\sum_{l=1}^{N_{00}} A^{N,\PP}_{ij}(T)A^{N,\PP}_{ik}(T)\left|f^{N,\CC}_{jl}(s)f^{N,\CC}_{kl}(t)\cov[b_l(s),b_l(t)]\right|\\
        &\quad+\sum_{j\in\calp^N} (A^{N,\PP}_{ij}(T))^2 \left|f^{N,\CC}_{j(N_{00}+i)}(s)f^{N,\CC}_{j(N_{00}+i)}(t)\cov[b_{N_{00}+i}(s),b_{N_{00}+i}(t)]\right| \Bigg)^{1/2}  \\
				&\leq \phi(T)  v_b(T)v_f(T) \frac{\sqrt{N_{00}} p^N_f(T)+\sqrt{N}}{R^N_A} \to 0,\\
				r^N_{11}(T)&= \sup_{i\in\bbn}\sup_{t\in[0,T]} \Bigg(\sum_{j,k\in\calp^N}\sum_{l=1}^{N_{00}} A^{N,\PP}_{ij}(T)A^{N,\PP}_{ik}(T)\left|\bbe[\rho^{N,\CC}_{jl}(t)\rho^{N,\CC}_{kl}(t)]\right|\var[M_l(1)]\\
        &\quad+\sum_{j\in\calp^N} (A^{N,\PP}_{ij}(T))^2 \left|\bbe[(\rho^{N,\CC}_{j(N_{00}+i)}(t))^2]\right|\var[M_{N_{00}+i}(1)] \Bigg)^{1/2}\\
				&\leq \phi(T)  v_{\rho,M}(T) \frac{\sqrt{N_{00}} p^N_\rho(T) + \sqrt{N}}{R^N_A} \to 0,
				\end{align*}
				and similarly $r^N_{10}(T),r^N_{12}(T)\to0$ as $N\to\infty$.
			\item[(12b)] $A^{N,2}$ and $\Si^{N,2}$ are sparse in the sense of \eqref{A-sparse}. Then one can deduce $r^N_\iota(T)\to0$ for $\iota=9,10,11,12$ as before.
		\eenu

We conclude this subsection with two remarks.

\brem In the sparseness conditions \eqref{L-uncorr}--\eqref{pNrho} it is not essential that the majority of entries is exactly zero. As one can see from the definition of the rates \eqref{rates}, they depend continuously on the underlying matrix entries. It suffices therefore that the matrix entries are small enough in a large proportion. \halmos
\erem

\brem\label{gen3}
What can be said about Theorem~\ref{LLN} in the general case of nonlinear Lipschitz coefficients $a^N$ and $\si^N$, apart from the special case discussed in Remark~\ref{gen2}? In fact, a law of large numbers in the fashion of Theorem~\ref{LLN} can still be shown, but under more stringent conditions: namely we have to require condition (10b) above in addition, with $A^{N,\PP}$ and $\Si^{N,\PP}$ now containing the Lipschitz constants of the interaction terms. The reason is that (10b) suffices to make $r^N_\iota$ ($\iota\in\{3,4,7, \ldots ,12\}$) small. The remaining four rates are unrelated to $a^N$ and $\si^N$ and therefore not affected by their nonlinear structure. It is important to notice that conditions like (10a) and (12a) are no longer sufficient to make the corresponding rates small. The reason is that they are conditions of correlation type. Since correlation is a \emph{linear} measure of dependence, it is not surprising that these conditions are not suitable for the nonlinear case. We do not go into the details at this point. \halmos
\erem

\subsection{Networks arising from preferential attachment}\label{Sect33}
As demonstrated in the last subsection, the crucial criterion for the rates \eqref{rates} in Theorem~\ref{LLN} to vanish asymptotically with growing network size can be described as a combination of sparse interaction and sparse correlation among the particles. Condition \eqref{A-sparse} plays a distinguished role here: when valid, it implies that eight out of twelve rates in \eqref{rates} are small. Moreover, it is the key factor for a nonlinear generalization of Theorem~\ref{LLN} to hold or not; see Remark~\ref{gen3}. The aim of this subsection is therefore to find algorithms for the generation of the underlying networks such that the resulting interaction matrices satisfy \eqref{A-sparse}. We will assume that $a^{N,\PP}(t)=a^{N,\PP}$ and $\si^{N,\PP}(t)=\si^{N,\PP}$ are independent of $t\in\bbr_+$, such that also $A^{N,\PP}(t)$ and $\Si^{N,\PP}(t)$ as well as $p^N_{A,1}(t)$ and $p^N_{\Si}(t)$ (see \eqref{A-sparse} for their definitions) are independent of $t$. Furthermore, we only concentrate on $p^N_{A,1}$ as the analysis for $p^N_\Si$ is completely analogous.

We will base the creation of the IPS network on dynamical random graph mechanisms. Since we are mainly interested in heterogeneous graphs, we will investigate the \emph{preferential attachment} or \emph{scale-free} random graph \citep{Barabasi99}. There are many similar but different constructions of preferential attachment graphs; in the following, we rely on the construction of \citep{Bollobas03} for directed graphs. We remark that the random graphs to be constructed will be indexed by $N$, corresponding to a family of growing networks for the IPS. In particular, ``time'' in the random graph process must not be confused with the time $t$ in the IPS \eqref{SDE}; the correct view is rather that the IPS network has been built from the random graphs before time $t=0$, and, of course, independently of all random variables in \eqref{SDE}.

The preferential attachment algorithm works as follows: we start with $G(0)=(V,E(0))$, a given graph consisting of vertices $V=\bbn$ and edges $E(0)=\{e_1, \ldots, e_\nu\}$, where $\nu \in\bbn$ and $e_i$ stands for a directed edge between two vertices. We allow for multiple edges and loops in our graphs. Without loss of generality, we assume that the set of vertices in $G(0)$ with at least one neighbour  given by $\{1,\ldots,n(0)\}$ with some $n(0)\in\bbn$. Furthermore, we fix $\al,\beta,\ga\in\bbr_+$ with $\al+\beta+\ga=1$ and $\al+\ga>0$ and two numbers $\delta^\inn,\delta^\out\in\bbr_+$. For $N\in\bbn$ we construct $G(N)=(V,E(N))$ from $G(N-1)$ according to the following algorithm.
\bit
  \item With probability $\al$, we create a new edge $e_{\nu+N}$ from $v=n(N-1)+1$ to a node $w$ that is already connected in $G(N-1)$. Here $w$ is chosen randomly from $\{1,\ldots, n(N-1)\}$ according to the probability mass function
  \[ \frac{\dd^\inn_{G(N-1)}(w)+\delta^\inn}{\nu+N-1+\delta^\inn n(N-1)},\quad w\in \{1,\ldots, n(N-1)\}, \]
  where $\dd^\inn_G(v)$ denotes the in-degree of vertex $v$ in a graph $G$. We define $n(N):=n(N-1)+1$ and $E(N):=E(N-1)\cup\{e_{\nu+N}\}$.
  \item With probability $\beta$, a new edge $e_{\nu+N}$ is formed from some vertex $v\in \{1,\ldots, n(N-1)\}$ to some $w\in \{1,\ldots, n(N-1)\}$ (the case $v=w$ is possible). Here $v$ and $w$ are chosen independently according to the probability mass functions
  \[ \frac{\dd^\out_{G(N-1)}(v)+\delta^\out}{\nu+N-1+\delta^\out n(N-1)},\quad \frac{\dd^\inn_{G(N-1)}(w)+\delta^\inn}{\nu+N-1+\delta^\inn n(N-1)},\quad v,w\in \{1,\ldots, n(N-1)\}, \]
  respectively, where $\dd^\out_G(v)$ denotes the out-degree of vertex $v$ in a graph $G$. Moreover, we set $n(N):=n(N-1)$ and $E(N):=E(N-1)\cup\{e_{\nu+N}\}$.
  \item With probability $\ga$, a new edge $e_{\nu+N}$ from some $v\in \{1,\ldots, n(N-1)\}$ to $w=n(N-1)+1$ is formed. Here $v$ is chosen randomly according to the probability mass function
  \[ \frac{\dd^\out_{G(N-1)}(v)+\delta^\out}{\nu+N-1+\delta^\out n(N-1)},\quad v\in \{1,\ldots, n(N-1)\}. \]
  We set $n(N):=n(N-1)+1$ and $E(N):=E(N-1)\cup\{e_{\nu+N}\}$.
\eit
Evidently, we always have $|E(N)|=\nu+N$ while the number $n(N)$ of non-isolated vertices in $G(N)$ is random in general.

The most important result for our purposes is the following one. We define
\[ M^\inn(N):=\max\{d^\inn_{G(N)}(i)\colon i\in\bbn\},\quad M^\out(N):=\max\{d^\out_{G(N)}(i)\colon i\in\bbn\},\quad N\in\bbn_0, \]
as the maximal in-degree and out-degree in $G(N)$, respectively.
\blem\label{prefatt} The maximum in-degree $M^\inn(N)$ and out-degree $M^\out(N)$ of $G(N)$ satisfy the following asymptotics:
\beq\label{asympmaxdeg} c^\inn(N)M^\inn(N) \to \mu^\inn,\quad c^\out(N)M^\inn(N) \to \mu^\out,\quad N\to\infty. \eeq
Here the convergence to the random variables $\mu^\inn$ and $\mu^\out$, respectively, holds in the almost sure as well as in the $L^p$-sense for all $p\in[1,\infty)$, and $(c^\inn(N))_{N\in\bbn}$ and $(c^\out(N))_{N\in\bbn}$ are sequences of random variables which can be chosen such that for every $\eps\in(0,\al+\ga)$ we have a.s.
\beq\label{c-asymp} c^\inn(N)^{-1} = \Omicron\left(N^\frac{\al+\beta}{1+\delta^\inn(\al+\ga-\eps)}\right),\quad c^\out(N)^{-1} = \Omicron\left(N^\frac{\beta+\ga}{1+\delta^\out(\al+\ga-\eps)}\right),\quad N\to\infty.\eeq
\elem
It follows from this lemma that for every $\eps\in(0,\al+\ga)$ we have a.s.
\[ M^\inn(N)=\Omicron\left(N^\frac{\al+\beta}{1+\delta^\inn(\al+\ga-\eps)}\right),\quad M^\out(N) = \Omicron\left(N^\frac{\beta+\ga}{1+\delta^\out(\al+\ga-\eps)}\right),\quad N\to\infty. \]
In particular, if $G(N)$ is used to model the underlying network of $a^N$ (i.e. an edge in $G(N)$ from $i$ to $j$ is equivalent to $a^N_{ij}\neq0$), we have
\[ p^N_{A,1} \leq M^\out(N) = \Omicron\left(N^\frac{\beta+\ga}{1+\delta^\out(\al+\ga-\eps)}\right),\quad N\to\infty. \]
In other words, the first part of condition \eqref{A-sparse} holds as soon as $R^N_A$, as specified through \eqref{Aspecial} and \eqref{sqrtrate}, increases in $N$ at least with rate 
\beq\label{RNA} N^\frac{\beta+\ga}{1+\delta^\out(\al+\ga-\eps)}  \eeq
for some small $\eps$. For example, in the classical case of Example~\ref{classex} where $R^N_A=N$, this is always true except in the case $\al=\delta^\out=0$, where all edges start from one of the initial nodes with probability one. We conclude that in all non-trivial situations of the preferential attachment model, the resulting networks are sparse enough for the law of large numbers implied by Theorem~\ref{LLN} to be in force.

\section{Large deviations}\label{Sect4}

In Theorem~\ref{LLN} we have established bounds on the mean squared difference between the IPS \eqref{SDE} and the PMFS \eqref{meanSDE}. In Sections~\ref{Sect31}--\ref{Sect33} we have given examples of dynamical networks in which these bounds converge to $0$ as the network size increases. A natural question is now whether a large deviation principle holds as $N\to\infty$, which would then assure that the probability of $X^N$ deviating strongly from $\bar X^N$ decreases exponentially fast in $N$. In the classical case of homogeneous networks, \citep{Dawson87} is the first paper to prove a large deviation principle for the empirical measures of the processes \eqref{McKean}. For heterogeneous networks, however, the empirical measure might no longer be a good quantity to investigate: the weight of a particle now depends on which particle's perspective is chosen. A sequence of differently weighted empirical measures seems to be more appropriate, but then their analysis becomes considerably more involved. Therefore, in this paper we take a more direct approach and study the large deviation behaviour of the difference $X^N-\bar X^N$ itself. In order to do so, we have to put stronger assumptions on the coefficients than in the previous sections. These are as follows. 

\benu
	\item[(A1)] $X^N(0)$ is deterministic for each $N\in\bbn$.
	\item[(A2)] For all $N\in\bbn$ we have $\si^N=0$. All other coefficients $a^{N,\CC}$, $a^{N,\PP}$, $\rho^{N,\CC}$ and $\rho^{N,\PP}$ are constant in time. $\rho^{N,\CC}$ and $\rho^{N,\PP}$ (resp. $a^{N,\CC}$ and $a^{N,\PP}$) only have $\ga(N)$ (resp. $\Ga(N)$) non-zero columns, where $\ga(N)$ forms a sequence of natural numbers increasing to infinity and $\Ga(N)$ grows at most like $\exp(\ga(N))$.
	\item[(A3)] All numbers in \eqref{finnum}, which are indexed by $N$ now, are bounded independently of $N$.
	\item[(A4)] $(M_i\colon i\in\bbn)$ is a sequence of independent mean-zero L\'evy processes whose Brownian motion part has variance $c_i$ and whose L\'evy measure is $\nu_i$. Moreover, there exists a real-valued mean-zero L\'evy process $M_0$ that dominates $M_i$, that is, its characteristics $c_0$ and $\nu_0$ satisfy $c_i\leq c_0$ and $\nu_i(A)\leq \nu_0(A)$ for all $i\in\bbn$ and Borel sets $A\subseteq\bbr$, and that has finite exponential moments of all orders: $\bbe[\ee^{uM_0(1)}]<\infty$ for all $u\in\bbr_+$.
	\item[(A5)] Assume that $G^N(t,s):=\ga(N)\ee^{a^N t}a^{N,\PP}\ee^{a^{N,\CC}s}\rho^{N,\CC}$, $s,t\in[0,T]$, converges uniformly to a limit $G(t,s)\in\bbr^{\bbn\times\bbn}$:
	\[ \sup_{i,j\in\bbn} \sup_{s,t\in[0,T]} |G_{ij}(t,s)|<\infty,\quad \sup_{i,j\in\bbn}\sup_{s,t\in[0,T]} |G^N_{ij}(t,s)-G_{ij}(t,s)|\to 0,\quad N\to\infty. \]
	\item[(A6)] With $R^N(t):=\ga(N)\ee^{a^N t}\rho^{N,\PP}$, $t\in[0,T]$, there exists $R(t)\in\bbr^{\bbn\times\bbn}$ such that
	\[ \sup_{i,j\in\bbn} \sup_{t\in[0,T]} |R_{ij}(t)|<\infty,\quad \sup_{i,j\in\bbn} \sup_{t\in[0,T]} |R^N_{ij}(t)-R_{ij}(t)|\to0,\quad N\to\infty.\]
	\item[(A7)] The following two quantities are finite:
	\begin{align*}
	q_1&:=\limsup_{N\to\infty} q_1(N):= \limsup_{N\to\infty} \sup_{i\in\bbn} \ga(N) \sum_{j=1}^\infty \sum_{k\neq j} |a^{N,\PP}_{ij}a^{N,\CC}_{jk}|,\\
	q_2&:=\limsup_{N\to\infty} q_2(N):=\limsup_{N\to\infty} \sup_{i,k\in\bbn} \ga(N) \sum_{j=1}^\infty |a^{N,\PP}_{ij}\rho^{N,\CC}_{jk}|.
	\end{align*}
	\item[(A8)] Define for $m\in\bbn\cup\{0\}$
	\[ \Psi_m(u):=\frac{1}{2}c_m u^2 + \int_\bbr (\ee^{uz}-1-uz)\,\nu_m(\dd z),\quad u\in\bbr_+.\]
	We assume that the following holds for every $d\in\bbn$: denoting for $m\in\bbn$, $r\in[0,T]$ and $\theta\in M^d_T$
	\[
		H_m(\theta,r):= \int_r^T \int_s^T \sum_{i=1}^d G_{im}(t-s,s-r) \,\theta_i(\dd t)\,\dd s+\int_r^T \sum_{i=1}^d R_{im}(t-r)\,\theta_i(\dd t),
	\]
	the sequence $\left(\int_0^T \Psi_m(H_m(\theta,r))\,\dd r\right)_{m\in\bbn}$ is Ces\`aro summable, i.e. the following limit exists:
	\beq\label{Cesaro} \lim_{N\to\infty} \frac{1}{\ga(N)} \sum_{m=1}^{\ga(N)} \int_0^T \Psi_m(H_m(\theta,r)) \,\dd r.  \eeq
\eenu

\bthm\label{LD} 
Let $T\in\bbr_+$. Under (A1)--(A8), the sequence $(X^N-\bar X^N)_{N\in\bbn}$ satisfies a large deviation principle in $(D^\infty_T, J_1)$ with a good rate function $I\colon D^\infty_T \to [0,\infty]$, that is, for every $\al\in\bbr_+$ the set $\{x\in D^\infty_T\colon I(x)\leq \al\}$ is compact in $D^\infty_T$ (with respect to the $J_1$-topology), and for every $M\in\cald^\infty_T$ we have
\[ -\inf_{x\in \mathrm{int}\,M} I(x) \leq \liminf_{N\to\infty} \frac{1}{\ga(N)}\log \bbp[X^N-\bar X^N \in M] \leq \limsup_{N\to\infty} \frac{1}{\ga(N)}\log \bbp[X^N-\bar X^N \in M] \leq -\inf_{x\in \mathrm{cl}\,M} I(x), \]
where $\mathrm{int}\,M$ and $\mathrm{cl}\,M$ denote the interior and the closure of $M$ in $(D^\infty_T,J_1)$, respectively. Moreover, the rate function $I$ is convex, attains its minimum $0$ uniquely at the origin and is infinite for $x\notin AC^\infty_T$.
\ethm

\brem\label{LD-ex} 
\benu
\item We cannot drop the requirement $\si^N=0$ or condition (A4) in Theorem~\ref{LD}. If violated, the processes $X^N$ and $\bar X^N$ will typically not have exponential moments of all order, whose existence is essential for our proof below. This kind of problem does not arise when empirical measures are considered as in \citep{Dawson87, Leonard95} for the homogeneous case.
\item The Ces\`aro summability condition (A8) accounts for the possible inhomogeneity of the coefficients and the distribution of the noises. It holds in particular for the homogeneous case. Since a convergent series is Ces\`aro summable with the same limit, it also holds when we have asymptotic homogeneity (in the sense that the sequence inside the sum of \eqref{Cesaro} converges with $m\to\infty$). 
\item With $\ga(N)=N$ and assumptions (A2)--(A4) in force, McKean's example \eqref{McKean} or the model considered in \citep{Kley14} both satisfy the assumptions of the theorem. In McKean's case our large deviation principle follows from that of \citep{Dawson87} for the empirical measure by applying the contraction principle.
\eenu
\erem


\section{Proofs}\label{Sect5}

We start with some preparatory results that are needed for the proof of Theorem~\ref{LLN}.

\blem\label{varX} Under the assumptions of Theorem~\ref{LLN} we have
\[ \sup_{i\in\bbn} \left\|\bar X_i^\ast(T)\right\|_{L^2} \leq V(T), \]
where $V(T)$ is given in Theorem~\ref{LLN}.
\elem
\bpr 
It is a consequence of \eqref{meanSDE} and the Burkholder-Davis-Gundy inequality that for all $t\in[0,T]$ and $i\in\bbn$
\begin{align*}
\big\|(\bar X_i)^\ast(t)\big\|_{L^2} &\leq \big\|X_i(0)\big\|_{L^2} + \int_0^t \sum_{j=1}^\infty A_{ij}(T)\big\|(\bar X_j)^\ast(s)\big\|_{L^2}\,\dd s\\
&\quad+ 2\var[L_i(1)]\left(\int_0^t \left(\sum_{j=1}^\infty \Si_{ij}(T)\big\|(\bar X_j)^\ast(s)\big\|_{L^2}\right)^2\,\dd s \right)^{1/2}\\
&\quad+ \int_0^t \sum_{j=1}^\infty \big|f_{ij}(s)\big| \big\|b_j(s)\big\|_{L^2}\,\dd s + 2\left(\sum_{j,k=1}^\infty \int_0^t \bbe\big[\rho^\CC_{ij}(s)\rho^\CC_{ik}(s)c_{jk}(s)\big] \,\dd s\right)^{1/2}.
\end{align*}
Therefore, if we define $w(t):=\sup_{i\in\bbn} \|(\bar X_i)^\ast(t)\|_{L^2}$, we obtain
\begin{align*}
w(t)&\leq v_X + v_f(T)v_b(T)T + 2v_{\rho,M}(T)T^{1/2} + v_a(T) \int_0^t w(s)\,\dd s + 2v_L v_\si(T) \left(\int_0^t (w(s))^2\,\dd s\right)^{1/2}\\
&\leq v_X + v_f(T)v_b(T)T + 2v_{\rho,M}(T)T^{1/2} + (v_a(T)T^{1/2} + 2v_L v_\si(T))\left(\int_0^t (w(s))^2\,\dd s\right)^{1/2}.
\end{align*}
Now we square the last inequality, apply the basic estimate $(a+b)^2\leq 2(a^2+b^2)$ and use Gronwall's inequality to deduce our claim, namely that
\[ w(T) \leq \sqrt{2}\ee^{(v_a(T)T^{1/2} + 2v_L v_\si(T))^2 T} \left(v_X + v_f(T)v_b(T)T + 2v_{\rho,M}(T)T^{1/2}\right). \] \halmos
\epr


\blem\label{barXsol} Let $T\in\bbr_+$ and assume the finiteness of the numbers \eqref{finnum}. We fix some $j\in\bbn$ throughout this lemma and define for $t\in[0,T]$
\begin{align*} Y_j(t)&:=Y_j^{1}+Y_j^{2}(t) + Y_j^{3}(t) + Y_j^{4}(t) + Y_j^{5}(t)\\
&:=(X_j(0) -\bbe[X_j(0)]) + \sum_{k\neq j} \int_0^t a^\CC_{jk}(s)(\bar X_k(s)-\bbe[\bar X_k(s)])\,\dd s \\
&\quad+ \sum_{k=1}^\infty \int_0^t \big(\si^\CC_{jk}(s)\bar X_k(s-)+\si^\PP_{jk}(s)\bbe[\bar X_k(s)]\big) \,\dd L_j(s)\\
&\quad+ \sum_{k=1}^\infty \int_0^t f^\CC_{jk}(s)(b_k(s)-\bbe[b_k(s)])\,\dd s + \sum_{k=1}^\infty \int_0^t \rho^\CC_{jk}(s)\,\dd M_k(s). \end{align*}
Furthermore, introduce the integrals
\beq\label{itint} I^j_0[x](t) := x(t),\quad I^j_{n}[x](t):=\int_0^t a^\CC_{jj}(s) I^j_{n-1}[x](s) \,\dd s,\quad n\in\bbn, \ \eeq
where $x\colon[0,T]\to\bbr$ is a measurable function such that the integrals in \eqref{itint} exist for $t\in[0,T]$. Then
\beq\label{Xjrep} \bar X_j(t) - \bbe[\bar X_j(t)] = \sum_{n=0}^\infty I^j_n[Y_j](t) = \sum_{\iota=1}^5 \sum_{n=0}^\infty I^j_n[Y^\iota_j](t),\quad t\in[0,T],  \eeq
where the sums converge with respect to the maximal $L^2$-norm $X\mapsto \|X^\ast(T)\|_{L^2}$.
\elem
\bpr We deduce from \eqref{meanSDE} that
\begin{align*}
&~\bar X_j(t)-\bbe[\bar X_j(t)] \\
=&~ (X_j(0)-\bbe[X_j(0)]) + \sum_{k\neq j} \int_0^t a^\CC_{jk}(s)(\bar X_k(s)-\bbe[\bar X_k(s)])\,\dd s \\
&~ + \int_0^t a^\CC_{jj}(s) (\bar X_j(s)-\bbe[\bar X_j(s)])\,\dd s + \sum_{k=1}^\infty \int_0^t \big(\si^\CC_{jk}(s)\bar X_k(s-)+\si^\PP_{jk}(s)\bbe[\bar X_k(s)]\big) \,\dd L_j(s) \\
&~ + \sum_{k=1}^\infty \int_0^t f^\CC_{jk}(s)(b_k(s)-\bbe[b_k(s)])\,\dd s + \sum_{k=1}^\infty \int_0^t \rho^\CC_{jk}(s)\,\dd M_k(s)\\
=&~ I^j_1[\bar X_j-\bbe[\bar X_j]](t) + Y_j(t).
\end{align*}
Iterating this equality $n$ times, we obtain
\beq\label{nit} \bar X_j(t)-\bbe[\bar X_j(t)] = \sum_{\nu=0}^{n-1} I^j_\nu[Y_j](t) + I^j_n[\bar X_j-\bbe[\bar X_j]](t),\quad t\in[0,T]. \eeq
Next, observe that for any c\`adl\`ag process $(X(t))_{t\in\bbr_+}$ with $\|X^\ast(T)\|_{L^2}<\infty$ we have
\[ \left\|(I^j_\nu[X])^\ast(T)\right\|_{L^2} \leq \|X^\ast(T)\|_{L^2} I^j_\nu[1](T) \leq \|X^\ast(T)\|_{L^2} \frac{(A^\CC_{jj}(T))^\nu}{\nu!}, \]
which is summable in $\nu$. Thus, recalling from Lemma~\ref{varX} that both $Y_j$ and $\bar X_j-\bbe[\bar X_j]$ have finite maximal $L^2$-norm, we can let $n\to\infty$ in \eqref{nit} and get
\[ \bar X_j(t)-\bbe[\bar X_j(t)] = \sum_{\nu=0}^\infty I^j_\nu[Y_j](t),\quad t\in[0,T], \]
which is the first assertion. The second part of formula \eqref{Xjrep} holds by linearity. \halmos
\epr

\noindent\bff{Proof of Theorem~\ref{LLN}}.\quad The existence and uniqueness of solutions to \eqref{SDE} and \eqref{meanSDE} follow from the general theory of SDEs, see \citep{Protter05}, Theorem~V.7. Since the numbers \eqref{finnum} are finite, there are no difficulties in dealing with infinite-dimensional systems as in our case.

It follows from \eqref{SDE} and \eqref{meanSDE} that the difference between $X$ and $\bar X$ satisfies the SDE
\begin{align*}
	\dd(X(t)-\bar X(t)) &= \left(a(t)(X(t)-\bar X(t)) + a^\PP(t)(\bar X(t) - \bbe[\bar X(t)])\right)\,\dd t \\
	&\quad+ \left(\si(t)(X(t-)-\bar X(t-)) + \si^\PP(t)(\bar X(t-) - \bbe[\bar X(t)])\right).\dd L(t) \\
	&\quad+ f^\PP(t)(b(t)-\bbe[b(t)])\,\dd t + \rho^\PP(t)\,\dd M(t),\quad t\in\bbr_+,\\
	X(0)-\bar X(0) &= 0.
\end{align*}
Thus, denoting the left-hand side of \eqref{LLN1} by $\Delta(T)$, we obtain from the Burkholder-Davis-Gundy inequality and Jensen's inequality that
\begin{align} \Delta(T) &\leq v_a(T) \int_0^T \Delta(t)\,\dd t + 2v_\si(T)v_L \left(\int_0^T (\Delta(t))^2\,\dd t\right)^{1/2}\nonumber\\
&\quad+\left|\int_0^T \left\|a^\PP(t)(\bar X(t)-\bbe[\bar X(t)])\right\|_{L^2}\,\dd t\right|_\infty + \bigg| \left\|\left(\si^\PP(\bar X-\bbe[\bar X])\,.\, L\right)^\ast(T) \right\|_{L^2} \bigg|_\infty\nonumber \\
&\quad+T \sup_{t\in[0,T]} \Big| \left\|f^\PP(t)(b(t)-\bbe[b(t)])\right\|_{L^2} \Big|_\infty + \Big| \left\|(\rho^\PP\cdot M)^\ast(T)\right\|_{L^2} \Big|_\infty \nonumber \\
&\leq (T^{1/2}v_a(T)+2v_\si(T)v_L)\left(\int_0^T (\Delta(t))^2\,\dd t\right)^{1/2} + \sum_{\iota=1}^4\Delta^\iota(T), \label{DeltaN}\end{align}
where $\Delta^\iota(T)$ stands for the last four summands in the line before. So Gronwall's inequality produces the bound
\beq \Delta(T) \leq K(T)\sum_{\iota=1}^4\Delta^\iota(T),\eeq
where $K(T)=\sqrt{2}\exp((T^{1/2}v_a(T)+2v_\si(T)v_L)^2 T)$. We now consider each $\Delta^\iota(T)$ separately. 

For $\iota=3$ we simply have
\beq\label{iota3} \Delta^3(T)\leq T \sup_{t\in[0,T]} \sup_{i\in\bbn} \left(\sum_{j,k=1}^\infty f^\PP_{ij}(t)f^\PP_{ik}(t)\cov[b_j(t),b_k(t)]\right)^{1/2} = T r_5(T).\eeq

For $\iota=4$ another application of the Burkholder-Davis-Gundy inequality yields
\begin{align} \Delta^4(T) &\leq 2\sup_{i\in\bbn} \left(\sum_{j,k=1}^\infty \bbe\left[\int_0^T \rho^\PP_{ij}(t)\rho^\PP_{ik}(t)\,\dd[M_j,M_k](t)\right]\right)^{1/2} \nonumber \\
&\leq  2T^{1/2}\sup_{i\in\bbn} \sup_{t\in[0,T]} \left(\sum_{j,k=1}^\infty \bbe\left[\rho^\PP_{ij}(t)c_{jk}(t)\rho^\PP_{ik}(t)\right]\right)^{1/2}=2T^{1/2} r_6(T). \label{2term}
\end{align}

For $\iota=1$, we use Lemma~\ref{barXsol} including the notations introduced there and the fact that for all stochastic processes $(X(t))_{t\in\bbr_+}$ and $(Y(t))_{t\in\bbr_+}$ with c\`adl\`ag sample paths we have
\beq \sup_{r,s\in[0,T]} \left| \bbe\left[I^j_n[X](s)I^k_m[Y](r)\right] \right| \leq \frac{(A^\CC_{jj}(T))^n}{n!} \frac{(A^\CC_{kk}(T))^m}{m!} \sup_{r,s\in[0,T]} |\bbe[X(s)Y(r)]|\label{rule2} \eeq
for any $j,k\in\bbn$ and $m,n\in\bbn\cup\{0\}$. In this way we obtain
\begin{align} \Delta^1(T)&=\left|\int_0^T \left\|a^\PP(t)(\bar X(t)-\bbe[\bar X(t)])\right\|_{L^2}\,\dd t\right|_\infty =\sup_{i\in\bbn} \int_0^T \left\|\sum_{j=1}^\infty a^\PP_{ij}(t)(\bar X_j(t)-\bbe[\bar X_j(t)])\right\|_{L^2}\,\dd t\nonumber\\
&\leq \sum_{\iota=1}^5 \sup_{i\in\bbn}\int_0^T \left\|\sum_{j=1}^\infty a^\PP_{ij}(t)\sum_{n=0}^\infty I^j_n[Y^\iota_j](t) \right\|_{L^2}\,\dd t\nonumber\\
&= \sum_{\iota=1}^5 \sup_{i\in\bbn} \int_0^T\left(\sum_{j,k=1}^\infty \sum_{n,m=0}^\infty a^\PP_{ij}(t) a^\PP_{ik}(t) \bbe\left[I^j_n[Y^\iota_j](t) I^k_m[Y^\iota_k](t)\right] \right)^{1/2}\,\dd t\nonumber\\
&\leq \ee^{|A^\CC(T)|_{\dd}} \sum_{\iota=1}^5 \sup_{i\in\bbn} \int_0^T \left(\sum_{j,k=1}^\infty A^\PP_{ij}(T)A^\PP_{ik}(T) \sup_{r,s\in[0,t]} \left|\bbe[Y^\iota_j(s)Y^\iota_k(r)]\right|\right)^{1/2}\,\dd t\nonumber\\
&=:\ee^{|A^\CC(T)|_{\dd}}\sum_{\iota=1}^5 R_{\iota}(T). \label{fiveterms}
\end{align}
Using Lemma~\ref{varX}, the five terms in \eqref{fiveterms} can be estimated as follows:
\begin{align*}
R_1(T)&\leq T \sup_{i\in\bbn} \left(\sum_{j,k=1}^\infty A^\PP_{ij}(T)A^\PP_{ik}(T)|\cov[X_j(0),X_k(0)]|\right)^{1/2} = Tr_1(T),\\
R_2(T)&\leq \sup_{i\in\bbn} \int_0^T \sum_{j=1}^\infty A^\PP_{ij}(T) \sup_{s\in[0,t]}\|Y^2_j(s)\|_{L^2} \,\dd t\\
&\leq \sup_{i\in\bbn} \int_0^T \sum_{j=1}^\infty A^\PP_{ij}(T) \left(\sum_{k\neq j} \int_0^t A^\CC_{jk}(s)\| \bar X_k(s)-\bbe[\bar X_k(s)]\|_{L^2}\,\dd s\right) \,\dd t\\
&\leq \frac{T^2}{2} V(T) \sup_{i\in\bbn} \sum_{j=1}^\infty \sum_{k\neq j} A^\PP_{ij} A^\CC_{jk} = \frac{T^2}{2}V(T)r_7(T), \\
R_3(T)&\leq \sup_{i\in\bbn} \int_0^T \Bigg(\sum_{j,k=1}^\infty A^\PP_{ij}(T)A^\PP_{ik}(T) \sup_{s\in[0,t]} \bigg|\bbe\Big[\Big(\left(\si^\CC \bar X+\si^\PP \bbe[\bar X]\right)_j\cdot L_j\Big)(s)\\
&\quad \times\Big(\left(\si^\CC\bar X+\si^\PP\bbe[\bar X]\right)_k\cdot L_k\Big)(s)\Big]\bigg| \Bigg)^{1/2}\,\dd t\\
&\leq \sup_{i\in\bbn} \int_0^T \Bigg(\sum_{j,k=1}^\infty A^\PP_{ij}(T)A^\PP_{ik}(T) \cov[L_j(1),L_k(1)]\\
&\quad\times\int_0^t \bbe\left[\left|\left(\si^\CC(s) \bar X(s)+\si^\PP(s) \bbe[\bar X(s)]\right)_j \left(\si^\CC(s)\bar X(s)+\si^\PP(s)\bbe[\bar X(s)]\right)_k \right|\right]\,\dd s \Bigg)^{1/2}\,\dd t\\
&\leq \frac{2}{3}T^{3/2}v_\si(T)V(T)r_3(T),\\
R_4(T)&\leq \sup_{i\in\bbn}\int_0^T \left(\sum_{j,k=1}^\infty A^\PP_{ij}(T)A^\PP_{ik}(T)\int_0^t \int_0^t \left|\sum_{l,m=1}^\infty f^\CC_{jl}(s) f^\CC_{km}(r) \cov[b_l(s),b_m(r)]\right|\,\dd r \,\dd s  \right)^{1/2} \,\dd t\\
&\leq \frac{T^2}{2} \sup_{i\in\bbn}\left(\sum_{j,k=1}^\infty A^\PP_{ij}(T)A^\PP_{ik}(T) \sup_{s,t\in[0,T]} \left|\sum_{l,m=1}^\infty f^\CC_{jl}(s) f^\CC_{km}(t) \cov[b_l(s),b_m(t)]\right|\right)^{1/2}\\
&=\frac{T^2}{2} r_9(T),\\
R_5(T)&\leq \sup_{i\in\bbn}\int_0^T\left(\sum_{j,k=1}^\infty A^\PP_{ij}(T)A^\PP_{ik}(T) \sup_{s\in[0,t]}\left|\sum_{l,m=1}^\infty \bbe\left[(\rho^\CC_{jl}\cdot M_l)(s)(\rho^\CC_{km}\cdot M_m)(s)\right]\right| \right)^{1/2}\,\dd t\\
&\leq \sup_{i\in\bbn}\int_0^T\left(\sum_{j,k=1}^\infty A^\PP_{ij}(T)A^\PP_{ik}(T) \int_0^t \left|\sum_{l,m=1}^\infty \bbe\left[\rho^\CC_{jl}(s) c_{lm}(s) \rho^\CC_{km}(s)\right]\right|\,\dd s \right)^{1/2}\,\dd t\\
&\leq\frac{2}{3}T^{3/2} r_{11}(T).
\end{align*}

The last step in the proof is the estimation of $\Delta^2(T)$. To this end, we make use of the Burkholder-Davis-Gundy inequality another time and get
\begin{align*}
\Delta^2(T)&\leq \sup_{i\in\bbn} \left\|\left(\big(\si^\PP(\bar X-\bbe[\bar X])\big)_i\cdot L_i\right)^\ast(T) \right\|_{L^2}\\
&\leq 2v_L \sup_{i\in\bbn} \left(\int_0^T  \bbe\left[ \big(\si^\PP(t)(\bar X(t)-\bbe[\bar X(t)])\big)_i^2\right]\,\dd t\right)^{1/2}.
\end{align*}
The further procedure is analogous to what we have done for $\Delta^1(T)$: instead of $a^\PP$ we have $\si^\PP$ here. We leave the details to the reader and only state the result, which is
\[ \Delta^2(T) \leq \sum_{\iota\in\{2,4,8,10,12\}} K_\iota(T)r_\iota(T). \]
This completes the proof of Theorem~\ref{LLN}.
\halmos

\vspace{\baselineskip}
Our next goal is to prove Lemma~\ref{prefatt} concerning the rate of growth of the maximal degree in the preferential attachment random graph as described in Section~\ref{Sect33}. For the undirected version as in \citep{Barabasi99} the corresponding result goes back to \citep{Mori05}. Indeed, the proof there basically works for our case as well, but there are some steps that require different arguments. Thus, we decided to include the proof to our lemma.

\vspace{\baselineskip}
\noindent\bff{Proof of Lemma~\ref{prefatt}}.\quad
The statement is evidently true for $M^\inn$ when $\al+\beta=0$ (resp. for $M^\out$ when $\beta+\ga=0$). In fact, for this extremal case, in every step of the random graph a new edge is created pointing to (resp. from) a new node. This means that $M^\inn(N)$ (resp. $M^\out(N)$) remains constant for all $N\in\bbn_0$, and the claim follows with $c^\inn=1$ (resp. $c^\out=1$) identically. In the other cases, we closely follow the proof of Theorem~3.1 in \citep{Mori05}. In addition to the notation introduced in Section~\ref{Sect33}, we further define for $N\in\bbn_0$ and $\diamond\in\{\inn,\out\}$:
\begin{align*}
S^\diamond(N)&:=\nu+N+\delta^\diamond n(N),\\
X^\diamond(N,j)&:=\dd^\diamond_{G(N)}(j)+\delta^\diamond,\quad j\in\bbn,\\
N^\diamond_j&:=\inf\{N\in\bbn_0\colon \dd^\diamond_{G(N)}(j)\neq0\},\quad j\in\bbn,\\
s^\diamond&:=\al\bone_{\{\diamond=\inn\}}+\beta+\ga\bone_{\{\diamond=\out\}},\\
c^\diamond(0,k)&:=1,\quad c^\diamond(N+1,k):=c^\diamond(N,k)\frac{S^\diamond(N)}{S^\diamond(N)+s^\diamond k},\quad k\in\bbr_+,\\
Z^\diamond(N,j,k)&:=c^\diamond(N^\diamond_j+N,k)\binom{X^\diamond(N^\diamond_j+N,j)+k-1}{k}\bone_{\{N^\diamond_j<\infty\}},\quad j\in\bbn,\quad k\in\bbr_+,\\
\calg(N)&:=\si(e_{\nu+i}\colon i=1,\ldots,N),\quad \calg(\infty):=\si\left(\bigcup_{N=0}^\infty \calg(N)\right).
\end{align*}
Obviously, $\calg(N)$ is the $\si$-field of all information up to step $N$ in the preferential attachment algorithm, and $N^\diamond_j$ is a stopping time relative to the filtration $(\calg(N))_{N\in\bbn}$ for every $j\in\bbn$. Analogously to Theorem~2.1 of \citep{Mori05} one can now show that for all $k\in\bbr_+$ and $j\in\bbn$ the sequence $(Z^\diamond(N,j,k))_{N\in\bbn_0}$ is a positive martingale relative to the filtration $(\calg(N^\diamond_j+N))_{N\in\bbn_0}$. As a consequence, Doob's martingale convergence theorem implies that
\beq\label{martconv} Z^\diamond(N,j,k) \to \zeta^\diamond(j,k)\quad\text{a.s.}\eeq
for some random variables $\zeta^\diamond(j,k)$. The convergence in \eqref{martconv} also holds in $L^p$ for all $p\in[1,\infty)$ because we have
\beq\label{ineqZ} Z^\diamond(N,j,k)^p \leq C(k,p) Z^\diamond(N,j,kp)\quad\text{a.s.} \eeq
for some deterministic constants $C(k,p)\in\bbr_+$ independent of $N$ and $j$. Indeed, on $\{N^\diamond_j<\infty\}$ we have by definition
\[ \frac{Z^\diamond(N,j,k)^p}{Z^\diamond(N,j,kp)} = \frac{c^\diamond(N^\diamond_j+N,k)^p}{c^\diamond(N^\diamond_j+N,kp)} {\binom{X^\diamond(N^\diamond_j+N,j) + k-1}{k}}^p{\binom{X^\diamond(N^\diamond_j+N,j) + kp -1}{kp}}^{-1}, \]
where
\begin{align*} &\frac{c^\diamond(N,k)^p}{c^\diamond(N,kp)} = \frac{c^\diamond(N-1,k)^p}{c^\diamond(N-1,kp)} \frac{S^\diamond(N)^{p-1} (S^\diamond(N)+s^\diamond kp)}{(S^\diamond(N)+s^\diamond k)^p}
\leq \frac{c^\diamond(N-1,k)^p}{c^\diamond(N-1,kp)} \leq \ldots \leq \frac{c^\diamond(0,k)^p}{c^\diamond(0,kp)} = 1,\\
&{\binom{x+k-1}{k}}^p{\binom{x+ kp -1}{kp}}^{-1}=\frac{\Ga(kp+1)}{\Ga(k+1)^p}\frac{\Ga(k+x)^p}{\Ga(x)^{p-1}\Ga(kp+x)} \stackrel{x\to\infty}{\longrightarrow} \frac{\Ga(kp+1)}{\Ga(k+1)^p},
\end{align*}
which shows \eqref{ineqZ}. Next, define for $N\in\bbn_0$ and $j\in\bbn$
\begin{align*} m^\diamond(N,j)&:=\max\{Z^\diamond(N-N^\diamond_i,i,1)\colon i=1,\ldots,j, N^\diamond_i\leq N\}, & m^\diamond(N)&:=m^\diamond(N,n(N)),\\
\mu^\diamond(j)&:=\max\{\zeta^\diamond(i,1)\colon i=1,\ldots,j\}, &\mu^\diamond&:=\sup\{\mu^\diamond(j)\colon j\in\bbn\},
\end{align*}
such that in particular the relationship $m^\diamond(N)=c^\diamond(N,1)(M^\diamond(N)+\delta^\diamond)$ holds. It is not hard to see that $(m^\diamond(N))_{N\in\bbn_0}$, as the maximum of martingale expressions, is a submartingale relative to $(\calg(N))_{N\in\bbn_0}$. By definition the sequence $(c^\diamond(N,k))_{N\in\bbn_0}$ decreases to $0$ as $N\to\infty$; more precisely, we have
\begin{align*}
c^\diamond(N,k)&=c^\diamond(N-1,k)\frac{S^\diamond(N-1)}{S^\diamond(N-1)+s^\diamond k}\leq c^\diamond(N-1,k) \frac{\nu+N-1+\delta^\diamond(n(0)+N-1)}{\nu+N-1+\delta^\diamond(n(0)+N-1)+s^\diamond k}\\
&\leq \prod_{j=0}^{N-1} \frac{(1+\delta^\diamond)j+\delta^\diamond n(0) + \nu}{(1+\delta^\diamond)j+\delta^\diamond n(0)+ \nu+s^\diamond k} = \frac{\Gamma\left(N+\frac{\delta^\diamond n(0)+\nu}{1+\delta^\diamond}\right)}{\Gamma\left(N+\frac{\delta^\diamond n(0)+\nu+s^\diamond k}{1+\delta^\diamond}\right)} \sim N^{-\frac{s^\diamond k}{1+\delta^\diamond}},\quad N\to\infty.
\end{align*}
As a consequence, when $p$ is large enough, 
\begin{align} \bbe[m^\diamond(N)^p] &\leq \bbe\left[\sum_{i=1}^{n(N)} Z^\diamond(N-N_i^\diamond,i,1)^p\right] \leq C(1,p) \bbe\left[\sum_{i=1}^{n(N)} Z^\diamond(N-N_i^\diamond,i,p)\right] \nonumber \\
&\leq C(1,p) \sum_{i=1}^\infty \bbe[Z^\diamond(0,i,p)] \leq C(1,p)\binom{n(0)+p+\delta^\diamond-1}{p}\sum_{i=1}^\infty \bbe[c^\diamond(N^\diamond_i,p)]\nonumber\\
&\leq C(1,p)\binom{n(0)+p+\delta^\diamond-1}{p}\left(n(0)+\sum_{i=1}^\infty \bbe[c^\diamond(i,p)]\right) <\infty \label{sumc}  \end{align}
independently of $N$. This implies that the submartingale $m^\diamond$ converges a.s. and in $L^p$ for all $p\in[1,\infty)$. It follows from \eqref{sumc} that for $j\geq n(0)$ we have
\begin{align} \bbe[(m^\diamond(N)-m^\diamond(N,j))^p] &\leq \bbe\left[\sum_{i=j+1}^{n(N)} Z^\diamond(N-N_i^\diamond,i,1)^p\right] \nonumber\\ 
&\leq C(1,p)\binom{n(0)+p+\delta^\diamond-1}{p}\sum_{i=j-n(0)+1}^\infty \bbe[c^\diamond(i,p)]. \label{diffm}
\end{align}
Letting $N\to\infty$, the left-hand side of \eqref{diffm} converges to 
\[\bbe\left[\left(\lim_{N\to\infty} c^\diamond(N,1) M^\diamond(N)-\mu^\diamond(j)\right)^p\right],\]
while the right-hand side is independent of $N$. Now taking the limit $j\to\infty$ and again assuming that $p$ is large, we obtain the desired result \eqref{asympmaxdeg}. Note at this point that $\mu^\diamond$ is indeed an a.s. finite random variable that belongs to $L^p$ for all $p\in[1,\infty)$, which is proved using a similar argument as in \eqref{sumc}. 

It remains to prove \eqref{c-asymp}. To this end, observe that by the law of large numbers we have $(n(N)-n(0))/N \to \al+\ga$ a.s. In other words, there exists for every $\eps\in(0,\al+\ga)$ a possibly random $\bar N\in\bbn$ such that for all $N\geq \bar N$ we have 
\[ \left|\frac{n(N)-n(0)}{N} - (\al+\ga)\right| \leq \eps, \]
or, equivalently, $n(N)\in[n(0)+(\al+\ga-\eps)N,n(0)+(\al+\ga+\eps)N]$. Consequently, for all $k\in\bbn$ and $N\geq \bar N$
\begin{align*}
c^\diamond(N,k)&=\prod_{i=0}^{N-1} \frac{S^\diamond(i)}{S^\diamond(i)+s^\diamond k} \geq \prod_{i=0}^{\bar N-1} \frac{S^\diamond(i)}{S^\diamond(i)+s^\diamond k} \prod_{i=\bar N}^{N-1} \frac{\nu+i+\delta^\diamond (n(0)+(\al+\ga-\eps)i)}{\nu+i+\delta^\diamond (n(0)+(\al+\ga-\eps)i)+s^\diamond k}\\
&=\prod_{i=0}^{\bar N-1} \frac{S^\diamond(i)}{S^\diamond(i)+s^\diamond k} \prod_{i=0}^{\bar N-1} \frac{\nu+i+\delta^\diamond (n(0)+(\al+\ga-\eps)i)+s^\diamond k}{\nu+i+\delta^\diamond (n(0)+(\al+\ga-\eps)i)}\\
&\quad\times\prod_{i=0}^{N-1} \frac{\nu+i+\delta^\diamond (n(0)+(\al+\ga-\eps)i)}{\nu+i+\delta^\diamond (n(0)+(\al+\ga-\eps)i)+s^\diamond k}\\
&=c^\diamond(\bar N,k) \frac{\Gamma\left(\bar N+\frac{\delta^\diamond n(0)+\nu+s^\diamond k}{1+\delta^\diamond(\al+\ga-\eps)}\right)\Gamma\left(N+\frac{\delta^\diamond n(0)+\nu}{1+\delta^\diamond(\al+\ga-\eps)}\right)}{\Gamma\left(\bar N+\frac{\delta^\diamond n(0)+\nu}{1+\delta^\diamond(\al+\ga-\eps)}\right)\Gamma\left(N+\frac{\delta^\diamond n(0)+\nu+s^\diamond k}{1+\delta^\diamond(\al+\ga-\eps)}\right)}\\
&\sim c^\diamond(\bar N,k) \frac{\Gamma\left(\bar N+\frac{\delta^\diamond n(0)+\nu+s^\diamond k}{1+\delta^\diamond(\al+\ga-\eps)}\right)}{\Gamma\left(\bar N+\frac{\delta^\diamond n(0)+\nu}{1+\delta^\diamond(\al+\ga-\eps)}\right)} N^{-\frac{s^\diamond k}{1+\delta^\diamond(\al+\ga-\eps)}},\quad N\to\infty.
\end{align*}
So choosing $c^\diamond(N):=c^\diamond(N,1)$ for $N\in\bbn_0$ fulfills \eqref{c-asymp}. \halmos

\vspace{\baselineskip} 
Finally, we turn to the proof of the large deviation result in Section~\ref{Sect4}. 

\vspace{\baselineskip}
\noindent\bff{Proof of Theorem~\ref{LD}}.\quad By definition, $X^N-\bar X^N$ satisfies the equation
\[ X^N(t)-\bar X^N(t) = \int_0^t a^N (X^N-\bar X^N)(s)\,\dd s + \int_0^t a^{N,\PP}(\bar X^N(s)-\bbe[\bar X^N(s)])\,\dd s + \rho^{N,\PP} M(t),\quad t\in[0,T],
\]
whose solution is
\[ X^N(t)-\bar X^N(t) = \int_0^t \ee^{a^N(t-s)} a^{N,\PP}(\bar X^N(s)-\bbe[\bar X^N(s)])\,\dd s + \int_0^t \ee^{a^N(t-s)}\rho^{N,\PP}\, \dd M(s),\quad t\in[0,T]. \]
In order to establish a large deviation principle, it suffices by Theorem~4.6.1 of \citep{Dembo10} to prove such a principle in $(D^d_T,\cald^d_T,J_1)$ for the first $d$ coordinates of the process for every $d\in\bbn$, that is, for the $D^d_T$-valued process
\begin{align}
 Y^N_i(t)&:=Y^{N,1}_i(t)+Y^{N,2}_i(t)+Y^{N,3}_i(t)\nonumber\\
 &:=\int_0^t \int_0^s \sum_{j,k=1}^\infty \ee^{a^N(t-s)}_{ij} a^{N,\PP}_{jk} \ee^{a^{N,\CC}_{kk}(s-r)} \sum_{l\neq k}  a^{N,\CC}_{kl} (\bar X^N_l(r)-\bbe[\bar X^N_l(r)])\,\dd r \,\dd s\nonumber\\
 &\quad+\int_0^t \sum_{j,k,l=1}^\infty \ee^{a^N(t-s)}_{ij}  a^{N,\PP}_{jk} \int_0^s \ee^{a^{N,\CC}_{kk}(s-r)}\rho^{N,\CC}_{kl}\,\dd M_l(r)\,\dd s \nonumber\\
 &\quad+\sum_{j,k=1}^\infty  \int_0^t \ee^{a^N(t-s)}_{ij} \rho^{N,\PP}_{jk} \,\dd M_k(s),\quad i=1,\ldots,d,\quad t\in[0,T], \label{YN} 
\end{align} 
where we have used the formula
\[
\bar X^N_i(t) - \bbe[\bar X^N_i(t)] =  \int_0^t \ee^{a^{N,\CC}_{ii}(t-s)}\sum_{j\neq i} a^{N,\CC}_{ij}(\bar X^N_j(s)-\bbe[\bar X^N_j(s)])\,\dd s + \sum_{j=1}^\infty \int_0^t \ee^{a^{N,\CC}_{ii}(t-s)} \rho^{N,\CC}_{ij}\,\dd M_j(s),
\]
valid for all $i\in\bbn$ and $t\in\bbr_+$. 
Actually, we will even prove the large deviation principle in $(D^d_T,\cald^d_T,U)$, which is stronger. To this end, we introduce the notation
\[ \hat x(t) := \sum_{k=1}^{[\ga(N) T]-1} x\left(\textstyle\frac{k}{\ga(N)}\right) \bone_{\left[\frac{k}{\ga(N)},\frac{k+1}{\ga(N)}\right)}(t) +  x\left(\textstyle\frac{[\ga(N) T]}{\ga(N)}\right) \bone_{\left[\frac{[\ga(N) T]}{\ga(N)},T\right)}(t),\quad t\in[0,T],\quad x\in D^d_T.\]
 Then by Theorem~4.2.13 of \citep{Dembo10} and Lemma~\ref{expeq} below we can equally well show a large deviation principle for $\hat Y^N=\hat Y^{N,1}+\hat Y^{N,2}+\hat Y^{N,3}$. The same principle will then hold for $Y^N$. But this is proved in Lemma~\ref{LD-d}. That the rate function for $X^N-\bar X^N$ is convex with unique minimum $0$ at $0$ and can only be finite for functions in $AC^\infty_T$, is inherited from the rate function of $\hat Y^N$.
\halmos

\blem\label{expeq} For each $d\in\bbn$ and $\iota=1, 2, 3$, the $D^d_T$-valued processes $Y^{N,\iota}$ and $\hat Y^{N,\iota}$ are \emph{exponentially equivalent}, that is, for all $\eps\in(0,1)$ we have
\[ \lim_{N\to\infty} \frac{1}{\ga(N)}\log \bbp\left[\sup_{t\in[0,T]} \sup_{i=1,\ldots,d} |Y^{N,\iota}_i(t)-\hat Y^{N,\iota}_i(t)| > \eps \right] = -\infty.\]
\elem
\bpr We start with $\iota=1$. Writing $\hat t = [\ga(N)t]/\ga(N)$ and $\diag(a):=a-a^\times$, we obtain
\begin{align} &~\sup_{t\in[0,T]}\sup_{i=1,\ldots,d} \left| Y^{N,1}_i(t)-\hat Y^{N,1}_i(t) \right| \nonumber\\
\leq&~ \sup_{t\in[0,T]}\left| (\ee^{a^N t}-\ee^{a^N \hat t})\left(\int_0^t \int_0^s \ee^{-a^N s}a^{N,\PP}\ee^{\diag(a^{N,\CC})(s-r)}(a^{N,\CC})^\times(\bar X^N(r)-\bbe[\bar X^N(r)])\,\dd r\,\dd s\right)\right|_\infty\nonumber\\
&+\sup_{t\in[0,T]} \left| \ee^{a^N \hat t}\left(\int_{\hat t}^t \int_0^s \ee^{-a^N s}a^{N,\PP}\ee^{\diag(a^{N,\CC})(s-r)}(a^{N,\CC})^\times(\bar X^N(r)-\bbe[\bar X^N(r)])\,\dd r\,\dd s\right)\right|_\infty. \label{split-two0}
\end{align}
We can proceed with these two terms separately. Since $|\ee^{a^N t}-\ee^{a^N \hat t}|_\infty \leq v_a \ee^{v_a T}/\ga(N)$, we have for the first term in
\eqref{split-two0}
\begin{align*}
&~\bbp\left[\sup_{t\in[0,T]}\left| (\ee^{a^N t}-\ee^{a^N \hat t})\left(\int_0^t \int_0^s \ee^{-a^N s}a^{N,\PP}\ee^{\diag(a^{N,\CC})(s-r)}(a^{N,\CC})^\times(\bar X^N(r)-\bbe[\bar X^N(r)])\,\dd r\,\dd s\right)\right|_\infty > \eps\right] \\
\leq&~ \bbp\left[\sup_{t\in[0,T]} \left|a^{N,\PP}(a^{N,\CC})^\times(\bar X^N(t)-\bbe[\bar X^N(t)]) \right|_\infty > \frac{\eps\ga(N)}{v_a(\ee^{v_{a}T})^3 T^2}\right] \\
\leq&~ \bbp\left[\sup_{t\in[0,T]}\sup_{i\in\bbn} \sum_{j=1}^\infty \sum_{k\neq j} \big|a^{N,\PP}_{ij}a^{N,\CC}_{jk}(\bar X^N_k(t)-\bbe[\bar X^N_k(t)])\big|  > \frac{\eps\ga(N)}{v_a(\ee^{v_{a}T})^3 T^2}\right]=:p(N).
\end{align*}
We note that $\xi^N:=\bar X^N-\bbe[\bar X^N]$ satisfies the integral equation
\[ \xi^N(t)=\int_0^t a^{N,\CC}\xi^N(s)\,\dd s + \rho^{N,\CC}M(t), \quad t\in\bbr_+.\]
Hence we have
\[ (\xi^N)^\ast(t) \leq \int_0^t |a^{N,\CC}| (\xi^N)^\ast(s) \,\dd s + (\rho^{N,\CC}M)^\ast(t),\quad t\in\bbr_+, \]
or after $n\in\bbn$ iterations,
\begin{align*} (\xi^N)^\ast(t) &\leq \frac{(|a^{N,\CC}|t)^n}{n!} (\xi^N)^\ast(t) + \sum_{m=0}^{n-1} \frac{(|a^{N,\CC}|t)^m}{m!}(\rho^{N,\CC}M)^\ast(t), \\
(\xi^N)^\ast(t)&\leq \left(\mathrm{I}-\frac{(|a^{N,\CC}|T)^n}{n!}\right)^{-1} \sum_{m=0}^{n-1} \frac{(|a^{N,\CC}|t)^m}{m!} (\rho^{N,\CC}M)^\ast(t),
\end{align*}
where the last line holds when $n$ is large enough such that $(v_a T)^n/n! < 1$. It is not difficult to recognize that the exact value of $n$ only affects some constants in the subsequent arguments with no impact on the final result; we therefore assume without loss of generality that $n=1$ (i.e. $v_a T<1$). Then 
\begin{align*}
p(N)&\leq \bbp\left[\sup_{i\in\bbn} \sum_{j=1}^\infty \sum_{k\neq j} \sum_{l=1}^{\infty} |a^{N,\PP}_{ij}a^{N,\CC}_{jk}| (\mathrm{I}-|a^{N,\CC}|T)^{-1}_{kl} \sup_{t\in[0,T]} \sum_{m=1}^{\ga(N)}|\rho^{N,\CC}_{lm}M_m(t)| >\frac{\eps\ga(N)}{v_a(\ee^{v_{a}T})^3 T^2}\right]\\
&\leq \bbp\left[ \frac{1}{\ga(N)}\sup_{l\in\bbn} \sup_{t\in[0,T]} \sum_{m=1}^{\ga(N)} |\rho^{N,\CC}_{lm} M_m(t)| >\frac{\eps\ga(N)(1-v_a T)}{q_1(N)v_a(\ee^{v_{a}T})^3 T^2}\right].
\end{align*}
Let $\la(N)$ be positive numbers to be chosen later. Using the independence of the L\'evy processes $M_i$ and Doob's maximal inequality, we arrive at
\begin{align*} p(N)&\leq \exp\left(-\frac{\eps\la(N)\ga(N)(1-v_a T)}{q_1(N)v_a(\ee^{v_{a}T})^3 T^2}\right) \prod_{m=1}^{\ga(N)}\bbe\left[\exp\left(\frac{\la(N)}{\ga(N)}  \sup_{l\in\bbn}|\rho^{N,\CC}_{lm}| |M_m(T)| \right)\right]\\
&\leq \exp\left(-\frac{\eps\la(N)\ga(N)(1-v_a T)}{q_1(N)v_a(\ee^{v_{a}T})^3 T^2}\right) \prod_{m=1}^{\ga(N)}\left(1+\bbe\left[\exp\left(\frac{\la(N)}{\ga(N)}  \sup_{l\in\bbn}|\rho^{N,\CC}_{lm}| M_m(T) \right)\right]\right)\\
&\leq 2^{\ga(N)}\exp\left(-\frac{\eps\la(N)\ga(N)(1-v_a T)}{q_1(N)v_a(\ee^{v_{a}T})^3 T^2}\right) \exp\left(T\ga(N) \Psi_0\left(\frac{\la(N)}{\ga(N)}\sup_{l,m\in\bbn} |\rho^{N,\CC}_{lm}|\right) \right).
\end{align*}
Now define 
\[ \la(N):=\ga(N)\Psi_0^{-1}(1)/\left(\sup_{l,m\in\bbn} |\rho^{N,\CC}_{lm}|\right),\quad N\in\bbn. \]
Since $\Psi_0$ is a convex function, its inverse $\Psi_0^{-1}$ is concave and therefore we have for large $N$ that $\la(N)\geq\Psi_0^{-1}(\ga(N))/\left(\sup_{l,m\in\bbn} |\rho^{N,\CC}_{lm}|\right)$, which increases to infinity with $N$. With this choice of $\la(N)$ it follows that
\[ \lim_{N\to\infty} \frac{1}{\ga(N)}\log p(N) = -\infty, \]
which completes the proof for the first term in \eqref{split-two0}. The second term can be treated in analogous way: now the factor $\ga(N)$ does not come from the difference $|\ee^{a^N t}-\ee^{a^N \hat t}|_\infty$, but from the domain of integration $(\hat t, t]$. The details are left to the reader.

For $\iota=2$ similar methods apply. Also here we do not give the details. Instead, we sketch the proof for $\iota=3$ where some modifications are necessary. Recalling the meaning of $\Ga(N)$ from (A2), we have
\begin{align} \sup_{t\in[0,T]}\sup_{i=1,\ldots,d} \left| Y^{N,3}_i(t)-\hat Y^{N,3}_i(t) \right|&\leq \sup_{t\in[0,T]} \sup_{i=1,\ldots,d} \left|\left( (\ee^{a^N t}-\ee^{a^N \hat t})\left(\int_0^t \ee^{-a^N s}\rho^{N,\PP}\,\dd M(s)\right)\right)_i\right| \nonumber\\
&\quad+\sup_{t\in[0,T]} \sup_{i=1,\ldots,d} \left|\left( \ee^{a^N \hat t}\left(\int_{\hat t}^t \ee^{-a^N s} \rho^{N,\PP}\,\dd M(s)\right) \right)_i \right|\nonumber\\
&\leq |\ee^{a^N t}-\ee^{a^N \hat t}|_\infty \sup_{t\in[0,T]}\sup_{i=1,\ldots,\Ga(N)} \left|\left( \int_0^t \ee^{-a^N s}\rho^{N,\PP}\,\dd M(s)\right)_i\right|\nonumber \\
&\quad+\ee^{v_{a}T}\sup_{t\in[0,T]}\sup_{i=1,\ldots,\Ga(N)} \left|\left( \int_{\hat t}^t \ee^{-a^N s} \rho^{N,\PP}\,\dd M(s)\right)_i\right|.\label{split-two}\end{align}
We can again consider these two terms separately. For the first one we have
\begin{align*} &~\bbp\left[|\ee^{a^N t}-\ee^{a^N \hat t}|_\infty \sup_{t\in[0,T]}\sup_{i=1,\ldots,\Ga(N)} \left|\left( \int_0^t \ee^{-a^N s}\rho^{N,\PP}\,\dd M(s)\right)_i\right|>\eps\right]\\
\leq&~ \Ga(N)\sup_{i=1,\ldots,\Ga(N)} \bbp\left[\sup_{t\in[0,T]} \left|\left(\int_0^t \ee^{-a^N s}\rho^{N,\PP}\,\dd M(s)\right)_i\right|>\frac{\eps\ga(N)}{v_a \ee^{v_a T}}\right]\\
\leq&~\Ga(N)\sup_{i=1,\ldots,\Ga(N)}\exp\left(-\frac{\eps\la(N)\ga(N)}{v_a \ee^{v_a T}}\right)\prod_{k=1}^{\ga(N)}\bbe\left[\exp\left(\la(N)\sup_{i\in\bbn} \left|\int_0^T \sum_{j=1}^\infty \ee^{-a^N s}_{ij}\rho^{N,\PP}_{jk}\,\dd M_k(s)\right|\right)\right]\\
=&~\Ga(N)\sup_{i=1,\ldots,\Ga(N)}\exp\left(-\frac{\eps\la(N)\ga(N)}{v_a \ee^{v_a T}}\right)\prod_{k=1}^{\ga(N)}\bbe\left[\exp\left(\frac{\la(N)}{\ga(N)}\sup_{i\in\bbn}\left|\int_0^T \sum_{j=1}^\infty \ee^{-a^N s}_{ij}\ga(N)\rho^{N,\PP}_{jk}\,\dd M_k(s)\right|\right)\right].
\end{align*}
Now recall that the stochastic integral in the last line has an infinitely divisible distribution. Moreover, the larger the integrand, the larger the exponential moment is. Since the integrand above is uniformly bounded in $i$ and $k$ by our hypotheses, the stochastic integral above can be replaced by some constant times $M_k(T)$ for the further estimation. Therefore, the remaining calculation can be completed as in the case $\iota=1$. For the second term in \eqref{split-two} the reasoning is the same, except that the factor $\ga(N)$ is now due to the domain $(\hat t, t]$ of the stochastic integral. Observe at this point that $M_k(t)- M_k(\hat t)$ has the same distribution as $M_k(t-\hat t)$ and that $|t-\hat t| \leq 1/\ga(N)$. Again, we do not carry out the details. \halmos
\epr

\blem\label{tight} For each $\iota=1,2,3$ the processes $(\hat Y^{N,\iota}\colon N\in\bbn)$ form an \emph{exponentially tight} sequence in $(D^d_T,\cald^d_T,U)$, that is, for every $L\in\bbr_+$ there exists a compact subset $K^L$ of $D^d_T$ (with respect to the uniform topology $U$) such that 
\[ \frac{1}{\ga(N)}\log\bbp[\hat Y^{N,\iota} \notin K^L] \leq -L. \]
\elem
\bpr
We first consider $\iota=1$. We will adapt the idea of Lemma~4.1 in \citep{deAcosta94} to our setting. As shown in part (I) of the proof there, it suffices to show that for every $a, \eps \in(0,\infty)$ there exist a compact set $H\subseteq D^d_T$, some $C\in(0,\infty)$ and $n\in\bbn$ such that for all $N\geq n$
\beq\label{equivtight} \bbp[ d(\hat Y^{N,1}, H) > \eps ] \leq C\ee^{-\ga(N)a}, \eeq
where $d(f,H):=\inf\{ \sup_{t\in[0,T]}\sup_{i=1,\ldots,d} |f_i(t)-g_i(t)|\colon g\in H \}$ for $f\in D^d_T$. In order to prove \eqref{equivtight}, we first define for $n\in\bbn$ and $A\subseteq\bbr^d$
\[ H_n(A):=\left\{f\in D^d_T\colon f=\sum_{\kappa=1}^{[\ga(n)T]-1} x_\kappa\bone_{\left[\frac{\kappa}{\ga(n)}, \frac{\kappa+1}{\ga(n)}\right)} + x_{[\ga(n)T]} \bone_{\left[\frac{\kappa}{\ga(n)}, T\right]},\quad x_1,\ldots, x_{[\ga(n)T]} \in A\right\}.  \]
It follows from Equation~(4.3) of \citep{deAcosta94} that for $N\geq n$, $A\subseteq \bbr^d$ and $f\in H_N(A)$ we have
\beq\label{dineq} d(f,H_n(A)) \leq \sup_{\kappa=0,\ldots,[\ga(n)T]-1} \sup_{\lambda\in\left[1, \frac{\ga(N)}{\ga(n)}+1\right)} \sup_{i=1,\ldots, d} \left|f_i\left(\frac{[\textstyle\frac{\ga(N)\kappa}{\ga(n)}]+\lambda}{\ga(N)}\wedge T\right)-f_i\left(\frac{[\textstyle\frac{\ga(N)\kappa}{\ga(n)}]}{\ga(N)}\right)\right|.  \eeq
Next, define $K:=[-1,1]^d$. Then for every $\beta\in(0,\infty)$ and $N\geq n$ we have
\beq \label{split} \bbp[d(\hat Y^{N,1},H_n(\beta K)) > \eps] \leq \bbp[\hat Y^{N,1}\notin H_N(\beta K)] + \bbp[\hat Y^{N,1} \in H_N(\beta K), d(\hat Y^{N,1}, H_n(\beta K))> \eps]. \eeq
The first probability is bounded as follows:
\begin{align*}
\bbp[\hat Y^{N,1}\notin H_N(\beta K)]&= \bbp\Bigg[\sup_{i=1,\ldots,d}\sup_{\kappa=1,\ldots,[\ga(N)T]} \Bigg|\int_0^{\frac{\kappa}{\ga(N)}}\int_0^s \sum_{j,k=1}^\infty \ee^{a^N(\kappa/\ga(N)-s)}_{ij}  a^{N,\PP}_{jk} \ee^{a^{N,\CC}_{kk}(s-r)}\\
&\quad\times \sum_{l\neq k} a^{N,\CC}_{kl} (\bar X^N_l(r)-\bbe[\bar X^N_l(r)])\,\dd r \,\dd s\Bigg| > \beta \Bigg]\\
&\leq  \bbp\left[\sup_{t\in[0,T]}\sup_{i=1,\ldots,d} \sum_{j=1}^\infty\sum_{k\neq j} \big|a^{N,\PP}_{ij} a^{N,\CC}_{jk} (\bar X^N_k(t)-\bbe[\bar X^N_k(t)])\big| > \frac{\beta}{(\ee^{v_a T}T)^2} \right]=:p^\prime(N).
\end{align*}
By the same arguments as in Lemma~\ref{expeq}, one obtains (again assuming $v_a T<1$ without loss of generality)
\begin{align*}
p^\prime(N)&\leq 2^{\ga(N)}\exp\left(-\frac{\beta\la(N)(1-v_a T)}{q_1(N)(\ee^{v_a T}T)^2}\right) \exp\left(T\ga(N) \Psi_0\left(\frac{\la(N)}{\ga(N)}\sup_{l,m\in\bbn} |\rho^{N,\CC}_{lm}|\right) \right).
\end{align*}
We choose $\la(N):=\ga(N)$ this time. Then we can make $\log(p^\prime(N))/\ga(N)$ arbitrarily small uniformly for large $N$ by varying the value of $\beta$. 

For the second step of the proof of \eqref{equivtight} we conclude from \eqref{dineq} that
\begin{align}
&~\bbp[\hat Y^{N,1} \in H_N(\beta K), d(\hat Y^{N,1}, H_n(\beta K))>\eps]\nonumber\\
\leq&~\bbp\Bigg[\sup_{\kappa=0,\ldots,[\ga(n)T]-1} \sup_{\lambda\in\left[1, \frac{\ga(N)}{\ga(n)}+1\right)} \sup_{i=1,\ldots, d} \left|\hat Y^{N,1}_i\left(\frac{[\textstyle\frac{\ga(N)\kappa}{\ga(n)}]+\lambda}{\ga(N)}\wedge T\right)-\hat Y^{N,1}_i\left(\frac{[\textstyle\frac{\ga(N)\kappa}{\ga(n)}]}{\ga(N)}\right)\right|>\eps\Bigg].\label{bigbracket}
\end{align}
For the further procedure, we split the difference in the last line into two terms in the same way as in \eqref{split-two0}. We only treat the corresponding first term. As before, the other one can be estimated similarly. Introducing the notation $t^{N,n}_{\kappa,\la}$ (resp. $t^{N,n}_{\kappa}$) for the time point in the first (resp. second) parenthesis of \eqref{bigbracket}, and observing that $0\leq \widehat{t^{N,n}_{\kappa,\la}}-t^{N,n}_{\kappa}\leq 2/\ga(N)+1/\ga(n)$, we obtain
\begin{align*}
&~\bbp\Bigg[\sup_{\kappa=0,\ldots,[\ga(n)T]-1}\sup_{\lambda\in\left[1, \frac{\ga(N)}{\ga(n)}+1\right)}\sup_{i=1,\ldots,d} \Bigg|\Bigg( \left(\ee^{a^N \widehat{t^{N,n}_{\kappa,\lambda}}}-\ee^{a^N t^{N,n}_\kappa}\right)\\
&\times\int_0^t\int_0^s \ee^{-a^N s}a^{N,\PP}\ee^{\diag(a^{N,\CC})(s-r)}(a^{N,\CC})^\times(\bar X^N(r)-\bbe[\bar X^N(r)])\,\dd r\,\dd s \Bigg)_i\Bigg|>\eps\Bigg]\\
\leq&~\bbp\left[\sup_{t\in[0,T]}\sup_{i=1,\ldots,d} \left| \sum_{j=1}^\infty \sum_{k\neq j} a^{N,\PP}_{ij} a^{N,\CC}_{jk}(\bar X^N_k(t)-\bbe[\bar X^N_k(t)]) \right|>\frac{\eps}{v_a(\ee^{v_{a}T})^3 T^2\left(\frac{2}{\ga(N)}+\frac{1}{\ga(n)}\right)}\right]\\
\leq&~ 2^{\ga(N)}\exp\left(-\frac{\eps\la(N)(1-v_a T)}{q_1(N)v_a(\ee^{v_{a}T})^3 T^2\left(\frac{2}{\ga(N)}+\frac{1}{\ga(n)}\right)}\right) \exp\left(T\ga(N) \Psi_0\left(\frac{\la(N)}{\ga(N)}\sup_{l,m\in\bbn} |\rho^{N,\CC}_{lm}|\right) \right), \end{align*}
where the last line follows in similar fashion as before. With $\la(N):=\ga(N)$ we can make, by taking $n$ large enough, the logarithm of the last term divided by $\ga(N)$ arbitrarily small for $N\geq n$. This finishes the proof for $\iota=1$. The case $\iota=2$ is analogous, while for $\iota=3$ the line of argument remains the same in principle, with slight changes to account for the discretization of L\'evy processes, cf. the proofs of Lemma~\ref{expeq} and Lemma~4.1 of \citep{deAcosta94}. \halmos
\epr

\blem\label{LD-d} The process $(\hat Y^N_i\colon i=1,\ldots,d)$ satisfies a large deviation principle in $(D^d_T, \cald^d_T, U)$ with a good convex rate function $I_d\colon D^d_T\to[0,\infty]$ such that $I_d(x)<\infty$ implies $x\in AC^d_T$. Moreover, we have $I_d(0)=0$ and this minimum is unique.
\elem
\bpr We apply the abstract G\"artner-Ellis theorem of \citep{deAcosta94}, Theorems~2.1 and 2.4, to $\hat Y^N$ and prove the following steps. 
\benu
	\item The laws of $\hat Y^N$, $N\in\bbn$, are exponentially tight in $(D^d_T,\cald^d_T,U)$.
	\item For all $\theta\in M^d_T$ 
	the limit $\La(\theta)=\lim_{N\to\infty} (\ga(N))^{-1}\La_N(\ga(N)\theta)$ exists, where
	\[ \La_N(\theta):=\log\bbe\left[\exp\left(\sum_{i=1}^d \int_0^T \hat Y^N_i(t)\,\theta_i(\dd t)\right)\right]. \] 
	\item The mapping $\La$ is $C^d_T$-G\^ateaux differentiable, in the sense that for all $\theta\in M^d_T$ there exists $x^\theta\in C^d_T$ such that for all $\eta\in M^d_T$
	\beq\label{Gateaux} \delta \La(\theta;\eta):=\lim_{\eps\to0} \frac{\La(\theta+\eps\eta)-\La(\theta)}{\eps} = \sum_{i=1}^d \int_0^T x^\theta_i(t)\,\eta_i(\dd t). \eeq
Part of the claim is that the limit in \eqref{Gateaux} exists. Moreover, we have $\La(0;\eta)=0$ for all $\eta\in M^d_T$.
	\item We have $\{x\in D^d_T\colon \La^\ast(x)<\infty\} \subseteq AC^d_T$, where 
	\[ \La^\ast(x):=\sup_{\theta\in M^d_T} \left(\sum_{i=1}^d \int_0^T x_i(t)\,\theta_i(\dd t) - \La(\theta)\right),\quad x\in D^d_T. \]
	\item For every $\al\in\bbr_+$ the set $\{x\in D^d_T\colon \La^\ast(x)\leq \al\}$ is compact in $(D^d_T,\cald^d_T,U)$.
\eenu
Part of the G\"artner-Ellis theorem is that the rate function $I_d$ is given by $\La^\ast$, the convex conjugate or Fenchel-Legendre transform of $\La$. Since $\La$ is a convex function in $\theta$ satisfying (3), the conjugate $\La^{\ast\ast}$ of $\La^\ast$ is again $\La$, see Theorem~12 of \citep{Rockafellar74}. Thus, by the first corollary to Theorem~1 in \citep{Asplund69}, we have $I_d(0) = \La^\ast(0)=0$ and this minimum is unique.

Let us now prove (1)--(5) above. Part (1) has been proved in Lemma~\ref{tight}. For (2) we first compute $\La$. For all $\theta\in M^d_T$ we have (recall that $\hat t:=[\ga(N)t]/\ga(N)]$)
\begin{align}
	\La_N(\ga(N)\theta)&= \log \bbe\Bigg[\exp\Bigg(\sum_{i=1}^d \int_0^T \Bigg(\int_0^{\hat t} \int_0^s \ga(N)\sum_{j,k,l,m=1}^\infty \ee^{a^N(\hat t-s)}_{ij} a^{N,\PP}_{jk}\ee^{a^{N,\CC}(s-r)}_{kl} \rho^{N,\CC}_{lm}\,\dd M_m(r)\,\dd s\nonumber\\
	&\quad+\ga(N)\sum_{j,k=1}^\infty  \int_0^{\hat t} \ee^{a^N(\hat t-s)}_{ij} \rho^{N,\PP}_{jk} \,\dd M_k(s)\Bigg)\,\theta_i(\dd t)\Bigg)\Bigg] \nonumber\\
	&=\sum_{m=1}^{\ga(N)}\log \bbe\Bigg[\exp\Bigg(\sum_{i=1}^d \int_0^T \Bigg(\int_0^{\hat t} \int_0^s \sum_{j,k,l=1}^\infty \ga(N)\ee^{a^N(\hat t-s)}_{ij}  a^{N,\PP}_{jk}\ee^{a^{N,\CC}(s-r)}_{kl} \rho^{N,\CC}_{lm}\,\dd M_m(r)\,\dd s\nonumber\\
	&\quad+  \int_0^{\hat t} \ga(N)\sum_{j=1}^\infty\ee^{a^N(\hat t-s)}_{ij} \rho^{N,\PP}_{jm} \,\dd M_m(s)\Bigg)\,\theta_i(\dd t) \Bigg)\Bigg] \label{La1}
\end{align}
by the independence of the processes $M_m$. By a stochastic Fubini argument (see Theorem~IV.65 of \citep{Protter05}), the term within the exponential in the previous line can also be written as ($\check s$ denotes the smallest multiple of $\ga(N)$ that is larger or equal to $s$)
\beq\label{HN1l} \int_0^{\hat T}\left(\int_r^{\hat T}\int_{\check s}^T \sum_{i=1}^d G^N_{im}(\hat t-s,s-r) \,\theta_i(\dd t)\,\dd s + \int_{\check r}^T \sum_{i=1}^d R^N_{im}(\hat t-r)\,\theta_i(\dd t) \right)\,M_m(\dd r), \eeq
and has an infinitely divisible distribution such that its logarithmic Laplace exponent in \eqref{La1} is explicitly known. Denoting the parenthesis in \eqref{HN1l} by $H^N_m(\theta,r)$, it is given by $\int_0^{\hat T} \Psi_m(H^N_m(\theta,r)) \,\dd r$. We claim that this term converges uniformly in $m$ to $\int_0^T \Psi_m(H_m(\theta,r)) \,\dd r$. Indeed, by the dominating property of $M_0$, the claim follows as soon as we can prove that $H^N_m(\theta,r)\to H_m(\theta,r)$ as $N\to\infty$, uniformly in $m\in\bbn$ and $r\in[0,T]$. This in turn follows from
\begin{align*}
&~|H_m(\theta,r)-H^N_m(\theta,r)|\\
\leq&~ \left| \int_{\hat T}^T \int_s^T \sum_{i=1}^d G_{im}(t-s, s-r) \,\theta_i(\dd t)\,\dd s\right| + \left|\int_r^{\hat T} \int_s^{\check s} \sum_{i=1}^d  G_{im}(t-s,s-r) \,\theta_i(\dd t)\,\dd s \right| \\
&+ \left|\int_r^{\hat T} \int_{\check s}^T \sum_{i=1}^d  (G_{im}(t-s,s-r)-G^N_{im}(t-s,s-r)) \,\theta_i(\dd t)\,\dd s\right|\\
&+ \left|\int_r^{\hat T} \int_{\check s}^T \sum_{i=1}^d  (G^N_{im}(t-s,s-r)-G^N_{im}(\hat t-s, s-r)) \,\theta_i(\dd t)\,\dd s\right|\\
&+ \left|\int_r^{\check r} \sum_{i=1}^d R_{im}(t-r)\,\theta_i(\dd t) \right| + \left| \int_{\check r}^T \sum_{i=1}^d (R_{im}(t-r)-R^N_{im}(t-r))\,\theta_i(\dd t) \right|\\
&+ \left|\int_{\check r}^T \sum_{i=1}^d (R^N_{im}(t-r)-R^N_{im}(\hat t-r))\,\theta_i(\dd t)\right|\\
\leq&~ d\sup_{i,m\in\bbn} \sup_{s,t\in[0,T]} |G_{im}(t,s)| \left(\ga(N)^{-1} \sup_{i=1,\ldots,d} |\theta_i|([0,T]) + \sup_{i=1,\ldots,d}\int_r^{\hat T} |\theta_i|([s,\check s)) \,\dd s\right) \\
&+d\sup_{i=1,\ldots,d} |\theta_i|([0,T])\Bigg(T\sup_{i,m\in\bbn} \sup_{s,t\in[0,T]} |G_{im}(t,s)-G^N_{im}(t,s)|+\frac{v_aT}{\ga(N)} \sup_{N,i,m\in\bbn} \sup_{s,t\in[0,T]} |G^N_{im}(t,s)|\\
&+\ga(N)^{-1}\sup_{i,j\in\bbn}\sup_{t\in[0,T]} |R_{ij}(t)| + \sup_{i,j\in\bbn}\sup_{t\in[0,T]} |R_{ij}(t)-R^N_{ij}(t)|+ \frac{v_a}{\ga(N)}\sup_{N,i,j\in\bbn}\sup_{t\in[0,T]} |R^N_{ij}(t)|\Bigg),
\end{align*}
where all terms converge to $0$ by hypothesis independently of $m$ and $r$. For the second summand one has to notice that the integral term equals $\int_r^{\hat T} \int_{\hat t}^t 1\,\dd s\,|\theta_i|(\dd t)$ and thus converges to $0$ uniformly in $i$ and $r$ with rate $1/\ga(N)$. Since the value of Ces\`aro sums remains unchanged under uniform approximations, it follows from assumption (A8) of Theorem~\ref{LD} that
\[ \La(\theta)=\lim_{N\to\infty} \frac{1}{\ga(N)}\sum_{m=1}^{\ga(N)} \int_0^T \Psi_m(H^N_m(\theta,r))\,\dd r=\lim_{N\to\infty} \frac{1}{\ga(N)}\sum_{m=1}^{\ga(N)} \int_0^T \Psi_m(H_m(\theta,r))\,\dd r.  \]

Next, we prove the $C^d_T$-G\^ateaux differentiability of $\La$. First, regarding the existence of $\delta \La(\theta;\eta)$ in \eqref{Gateaux}, we note that the mappings $M^d_T \to L^\infty(\bbn\times[0,T])$, $\theta\mapsto (H_m(\theta,r)\colon m\in\bbn,~r\in[0,T])$, are continuous linear operators and therefore Fr\'echet differentiable, which is stronger than G\^ateaux differentiability. Together with the fact that $\Psi_m$ is differentiable with locally bounded derivative, and $c_m\leq c_0$ and $\nu_m\leq \nu_0$ for all $m\in\bbn$, this implies that for every $\theta$ and $\eta$
\[ \eps^{-1}\int_0^T \left(\Psi_m(H_m(\theta+\eps\eta,r)) - \Psi_m(H_m(\theta,r))\right) \,\dd r \]
converges uniformly in $m\in\bbn$ as $\eps\to0$. This in turn proves the G\^ateaux differentiability of $\La$. Moreover, it enables us to compute the derivative explicitly. Using the chain rule for Fr\'echet derivatives, we obtain
\begin{align*}
&~\delta \La(\theta;\eta) = \lim_{N\to\infty}\sum_{m=1}^{\ga(N)} \int_0^T \lim_{\eps\to0} \frac{\Psi_m(H_m(\theta+\eps\eta,r))-\Psi_m(H_m(\theta,r))}{\eps}\,\dd r\\
=&~ \lim_{N\to\infty}\sum_{m=1}^{\ga(N)} \int_0^T \Psi^\prime_m(H_m(\theta,r))H_m(\eta,r)\,\dd r\\
=&~ \lim_{N\to\infty}\sum_{m=1}^{\ga(N)} \int_0^T \Psi^\prime_m(H_m(\theta,r))\left(\int_r^T \int_s^T \sum_{i=1}^d G_{im}(t-s,s-r) \,\eta_i(\dd t)\,\dd s +\int_r^T \sum_{i=1}^d R_{im}(t-r)\,\eta_i(\dd t) \right)\,\dd r\\
=&~ \lim_{N\to\infty}\sum_{m=1}^{\ga(N)} \sum_{i=1}^d \int_0^T \left(\int_0^t \int_0^s \Psi^\prime_m(H_m(\theta,r)) G_{im}(t-s,s-r) \,\dd r\,\dd s + \int_0^t  R_{im}(t-r) \,\dd r \right)\,\eta_i(\dd t)\\
=&~ \sum_{i=1}^d \int_0^T \int_0^t \lim_{N\to\infty}\sum_{m=1}^{\ga(N)}\left(\int_0^s  G_{im}(t-s,s-r)\Psi^\prime_m(H_m(\theta,r))\,\dd r +   R_{im}(t-s)\Psi^\prime_m(H_m(\theta,s)) \right)\,\dd s\,\eta_i(\dd t),
\end{align*}
where all interchanges of integration, summation and taking limits are justified by dominated convergence. From the last line we deduce the existence of $x^\theta\in C^d_T$ satisfying \eqref{Gateaux}. Since $H_m(0,r)=0$ and $\Psi^\prime_m(0)=0$, we have $\delta \La(0;\eta)=0$ identically.

Next, we demonstrate (4), namely that $\La^\ast$ only assumes finite values on the set $AC^d_T$, that is, $\La^\ast(x)<\infty$ implies that for every $\eps\in(0,\infty)$ there exists $\delta\in(0,\infty)$ such that $\sum_{i=1}^d \sum_{j=1}^n  |x_i(b_j)-x_i(a_j)| < \eps$ whenever $n\in\bbn$, $0\leq a_1<b_1 \leq \ldots \leq a_n < b_n\leq T$ and $\sum_{j=1}^n (b_j-a_j) < \delta$. In order to do so, we follow the strategy of proof in \citep{deAcosta94}, Theorem~3.1. We consider $\theta_i:=\sum_{j=1}^n \xi^j_i(\delta_{b_j}-\delta_{a_j})$ where $\xi^j_i\in\bbr^d$ is arbitrary. Then we evidently have $\theta_i((r,T])=\sum_{j=1}^n \xi^j_i \bone_{[a_j,b_j)}(r)$. Denoting $C_T:=T\sup_{s,t\in[0,T]} \sup_{i,m\in\bbn} |G_{im}(t,s)| + \sup_{t\in[0,T]} \sup_{i,m\in\bbn} |R_{im}(t)|$, it follows that
\begin{align*}
\La(\theta)&\leq \sup_{m\in\bbn} \int_0^T \Psi_m(H_m(\theta,r)\,\dd r \leq \int_0^T \Psi_0\left(C_T \sum_{i=1}^d \theta_i((r,T])\right)\,\dd r\\
& = \int_0^T \sum_{j=1}^n \Psi_0\left(C_T \sum_{i=1}^d \xi^j_i\right) \bone_{[a_j,b_j)}(r)\,\dd r \leq \sup_{j=1,\ldots,n}\Psi_0\left(C_T \sum_{i=1}^d |\xi^j_i|\right) \sum_{j=1}^n(b_j-a_j)\\
&=: C(T,\|\xi^1\|_1,\ldots,\|\xi^n\|_1) \sum_{j=1}^n(b_j-a_j), 
\end{align*}
where $\|\xi^j\|_1:=\sum_{i=1}^d |\xi^j_i|$.
As a consequence, we deduce from the definition of $\La^\ast$ that for all $\tau\in(0,\infty)$ and $\|\xi^j\|_1\leq \tau$
\[ \sum_{i=1}^d \sum_{j=1}^n \xi^j_i(x_i(b_j)-x_i(a_j)) \leq C(T,\tau,\ldots,\tau)\sum_{j=1}^n(b_j-a_j) + \La^\ast(x).  \]
Taking $\xi^j_i$ as the $\tau$ times the sign of $x_i(b_j)-x_i(a_j)$, it follows that
\beq\label{equicont} \sum_{i=1}^d \sum_{j=1}^n |x_i(b_j)-x_i(a_j)| \leq \tau^{-1}C(T,\tau,\ldots,\tau)\sum_{j=1}^n(b_j-a_j) + \tau^{-1}\La^\ast(x). \eeq
If $\La^\ast(x)<\infty$, we can now choose $\tau$ first and then $\delta$ to make the left-hand side arbitrarily small. 

It only remains to prove (5), the compactness of the level sets of $\La^\ast$. By step (4) and the lower semicontinuity of $\La^\ast$, its level sets are closed subsets of $AC^d_T$. Thus, the Arzel\`a-Ascoli theorem provides a compactness criterion. First, observe that for all $t\in[0,T]$ we have for $x\in AC^d_T$ with $\La^\ast(x)\leq \al$ that
\[ \sum_{i=1}^d |x_i(t)| = \sup_{\theta\in\Theta_t} \sum_{i=1}^d \int_0^T x_i(t)\,\theta_i(\dd t) \leq \al + \sup_{\theta\in\Theta_t} \La(\theta) < \infty, \]
where $\Theta_t$ is the finite collection of all $\theta$ for which each coordinate is either $\delta_t$ or $-\delta_t$. Second, for the proof of the uniform equicontinuity of the functions $x\in AC^d_T$ with $\La^\ast(x)\leq \al$, we recall from \eqref{equicont} that
\[ \sum_{i=1}^d |x_i(t)-x_i(s)| \leq \tau^{-1}C(T,\tau,\ldots,\tau)(t-s) + \tau^{-1}\al, \]
which converges to $0$ independently of $x$ when $s\uparrow t$ and $\tau\to\infty$. \halmos
\epr

\subsection*{Acknowledgement}
We thank Jean-Dominique Deuschel and Nina Gantert for their advice on the topic of large deviations. 
We are also grateful to Oliver Kley for drawing our attention to preferential attachment models. The first author further acknowledges support from the graduate program TopMath at Technische Universit\"at M\"unchen and the Studienstiftung des deutschen Volkes.

\addcontentsline{toc}{section}{References}
\bibliographystyle{plainnat}
\bibliography{bib-SR}

\begin{thebibliography}{37}
\providecommand{\natexlab}[1]{#1}
\providecommand{\url}[1]{\texttt{#1}}
\expandafter\ifx\csname urlstyle\endcsname\relax
  \providecommand{\doi}[1]{doi: #1}\else
  \providecommand{\doi}{doi: \begingroup \urlstyle{rm}\Url}\fi

\bibitem[Ahn and Ha(2010)]{Ahn10}
S.M. Ahn and S.-Y. Ha.
\newblock Stochastic flocking dynamics of the {C}ucker-{S}male model with
  multiplicative white noises.
\newblock \emph{J. Math. Phys.}, 51\penalty0 (10):\penalty0 103301, 2010.

\bibitem[Asplund and Rockafellar(1969)]{Asplund69}
E.~Asplund and R.T. Rockafellar.
\newblock Gradient of convex functions.
\newblock \emph{Trans. Amer. Math. Soc.}, 139:\penalty0 443--467, 1969.

\bibitem[Barab{\'a}si and Albert(1999)]{Barabasi99}
A.-L. Barab{\'a}si and R.~Albert.
\newblock Emergence of scaling in random networks.
\newblock \emph{Science}, 286\penalty0 (5439):\penalty0 509--512, 1999.

\bibitem[Battiston et~al.(2012)Battiston, Gatti, Gallegati, Greenwald, and
  Stiglitz]{Battiston12}
S.~Battiston, D.~Delli Gatti, M.~Gallegati, B.~Greenwald, and J.E. Stiglitz.
\newblock Liaisons dangereuses: Increasing connectivity, risk sharing, and
  systemic risk.
\newblock \emph{J. Econ. Dyn. Control}, 36\penalty0 (8):\penalty0 1121--1141,
  2012.

\bibitem[Bolley et~al.(2011)Bolley, Ca{\~n}izo, and Carrillo]{Bolley11}
F.~Bolley, J.A. Ca{\~n}izo, and J.A. Carrillo.
\newblock Stochastic mean-field limit: non-{L}ipschitz forces and swarming.
\newblock \emph{Math. Model. Methods Appl. Sci.}, 21\penalty0 (11):\penalty0
  2179--2210, 2011.

\bibitem[Bollob{\'a}s et~al.(2003)Bollob{\'a}s, Borgs, Chayes, and
  Riordan]{Bollobas03}
B.~Bollob{\'a}s, C.~Borgs, J.~Chayes, and O.~Riordan.
\newblock Directed scale-free graphs.
\newblock In \emph{Proceedings of the Fourteenth Annual ACM-SIAM Symposium on
  Discrete Algorithms (Baltimore, MD, 2003)}, pages 132--139, New York, 2003.
  ACM.

\bibitem[Bouchaud and M{\'e}zard(2000)]{Bouchaud00}
J.-P. Bouchaud and M.~M{\'e}zard.
\newblock Wealth condensation in a simple model of economy.
\newblock \emph{Phys. A}, 282\penalty0 (3--4):\penalty0 536--545, 2000.

\bibitem[Buckdahn et~al.(2014)Buckdahn, Li, and Peng]{Buckdahn14}
R.~Buckdahn, J.~Li, and S.~Peng.
\newblock Nonlinear stochastic differential games involving a major player and
  a large number of collectively acting minor agents.
\newblock \emph{SIAM J. Control Optim.}, 52\penalty0 (1):\penalty0 451--492,
  2014.

\bibitem[Budhiraja and Wu(2016)]{Budhiraja15}
A.~Budhiraja and R.~Wu.
\newblock Some fluctuation results for weakly interacting multi-type particle
  systems.
\newblock \emph{Stoch. Process. Appl.}, 126\penalty0 (8):\penalty0 2253--2296,
  2016.

\bibitem[Carmona and Zhu()]{Carmona14}
R.~Carmona and X.~Zhu.
\newblock A probabilistic approach to mean field games with major and minor
  players.
\newblock \emph{Ann. Appl. Probab.}, 26\penalty0 (3):\penalty0 1535--1580.

\bibitem[Dawson and G{\"a}rtner(1987)]{Dawson87}
D.A. Dawson and J.~G{\"a}rtner.
\newblock Large deviations from the {M}c{K}ean-{V}lasov limit for weakly
  interacting diffusions.
\newblock \emph{Stochastics}, 20\penalty0 (4):\penalty0 247--308, 1987.

\bibitem[de~Acosta(1994)]{deAcosta94}
A.~de~Acosta.
\newblock Large deviations for vector-valued {L}{\'e}vy processes.
\newblock \emph{Stoch. Process. Appl.}, 51\penalty0 (1):\penalty0 75--115,
  1994.

\bibitem[Degond et~al.(2014)Degond, Liu, and Ringhofer]{Degond14}
P.~Degond, J.-G. Liu, and C.~Ringhofer.
\newblock Evolution of the distribution of wealth in an economic environment
  driven by local {N}ash equilibria.
\newblock \emph{J. Stat. Phys.}, 154\penalty0 (3):\penalty0 751--780, 2014.

\bibitem[{Del Moral} et~al.(2001){Del Moral}, Kallel, and Rowe]{DelMoral01}
P.~{Del Moral}, L.~Kallel, and J.~Rowe.
\newblock Modelling genetics algorithms with interacting particle systems.
\newblock \emph{Rev. Mat. Teor. Apl.}, 8\penalty0 (2):\penalty0 19--77, 2001.

\bibitem[Dembo and Zeitouni(2010)]{Dembo10}
A.~Dembo and O.~Zeitouni.
\newblock \emph{Large Deviations Techniques and Applications}.
\newblock Springer, New York, 2nd edition, 2010.

\bibitem[Finnoff(1993)]{Finnoff93}
W.~Finnoff.
\newblock Law of large numbers for a general system of stochastic differential
  equations with global interaction.
\newblock \emph{Stoch. Process. Appl.}, 46\penalty0 (1):\penalty0 153--182,
  1993.

\bibitem[Finnoff(1994)]{Finnoff94}
W.~Finnoff.
\newblock Law of large numbers for a heterogeneous system of stochastic
  differential equations with strong local interaction and economic
  applications.
\newblock \emph{Ann. Appl. Probab.}, 4\penalty0 (2):\penalty0 494--528, 1994.

\bibitem[Fouque and Sun(2013)]{Fouque13}
J.-P. Fouque and L.-H. Sun.
\newblock Systemic risk illustrated.
\newblock In J.-P. Fouque and J.A. Langsam, editors, \emph{Handbook on Systemic
  Risk}, pages 444--452. Cambridge University Press, Cambridge, 2013.

\bibitem[Fournier and L{\"o}cherbach(2016)]{Fournier15}
N.~Fournier and E.~L{\"o}cherbach.
\newblock On a toy model of interacting neurons.
\newblock \emph{Ann. Inst. Henri Poincar{\'e} Probab. Stat.}, 52\penalty0
  (4):\penalty0 1844--1876, 2016.

\bibitem[G{\"a}rtner(1988)]{Gaertner88}
J.~G{\"a}rtner.
\newblock On the {M}c{K}ean-{V}lasov limit for interacting diffusions.
\newblock \emph{Math. Nachr.}, 137\penalty0 (1):\penalty0 197--248, 1988.

\bibitem[Giesecke et~al.(2015)Giesecke, Spiliopoulos, Sowers, and
  Sirignano]{Giesecke15}
K.~Giesecke, K.~Spiliopoulos, R.B. Sowers, and J.~Sirignano.
\newblock Large portfolio asymptotics for loss from default.
\newblock \emph{Math. Financ.}, 25\penalty0 (1):\penalty0 77--114, 2015.

\bibitem[Ichinomiya(2012)]{Ichinomiya12}
T.~Ichinomiya.
\newblock Bouchard-{M}{\'e}zard model on a random network.
\newblock \emph{Phys. Rev. E}, 86\penalty0 (3):\penalty0 036111, 2012.

\bibitem[Kley et~al.(2015)Kley, Kl{\"u}ppelberg, and Reichel]{Kley14}
O.~Kley, C.~Kl{\"u}ppelberg, and L.~Reichel.
\newblock Systemic risk through contagion in a core--periphery structured
  banking network.
\newblock In A.~Palczewski and {\L}.~Stettner, editors, \emph{Advances in
  Mathematics of Finance}, pages 133--149. Banach Center Publications,
  Warschau, 2015.

\bibitem[L{\'e}onard(1990)]{Leonard90}
C.~L{\'e}onard.
\newblock Some epidemic systems are long range interacting particle systems.
\newblock In J.-P. Gabriel, C.~Lef{\`e}vre, and P.~Picard, editors,
  \emph{Stochastic Processes in Epidemic Theory}, pages 170--183. Springer,
  Berlin, 1990.

\bibitem[L{\'e}onard(1995)]{Leonard95}
C.~L{\'e}onard.
\newblock Large deviations for long range interacting particle systems with
  jumps.
\newblock \emph{Ann. Inst. Henri Poincar{\'e} Probab. Stat.}, 31\penalty0
  (2):\penalty0 289--323, 1995.

\bibitem[{McKean, Jr.}(1966{\natexlab{a}})]{McKean66a}
H.P. {McKean, Jr.}
\newblock A class of {M}arkov processes associated with nonlinear parabolic
  equations.
\newblock \emph{Proc. Nat. Acad. Sci.}, 56\penalty0 (6):\penalty0 1907--1911,
  1966{\natexlab{a}}.

\bibitem[{McKean, Jr.}(1966{\natexlab{b}})]{McKean66b}
H.P. {McKean, Jr.}
\newblock Speed of approach to equilibrium for {K}ac's caricature of a
  {M}axwellian gas.
\newblock \emph{Arch. Ration. Mech. Anal.}, 21\penalty0 (5):\penalty0 343--367,
  1966{\natexlab{b}}.

\bibitem[{McKean, Jr.}(1967)]{McKean67}
H.P. {McKean, Jr.}
\newblock Propagation of chaos for a class of non-linear parabolic equations.
\newblock In \emph{Stochastic Differential Equations (Lecture Series in
  Differential Equations, Session 7, Catholic University, 1967)}, pages 41--57.
  Air Force Office of Scientific Research, Arlington, 1967.

\bibitem[M{\'o}ri(2005)]{Mori05}
T.F. M{\'o}ri.
\newblock The maximum degree of the {B}arab{\'a}si-{A}lbert random tree.
\newblock \emph{Comb. Probab. Comput.}, 14\penalty0 (3):\penalty0 339--348,
  2005.

\bibitem[Nagasawa and Tanaka(1987{\natexlab{a}})]{Nagasawa87}
M.~Nagasawa and H.~Tanaka.
\newblock Diffusion with interactions and collisions between coloured particles
  and the propagation of chaos.
\newblock \emph{Probab. Theory Relat. Fields}, 74\penalty0 (2):\penalty0
  161--198, 1987{\natexlab{a}}.

\bibitem[Nagasawa and Tanaka(1987{\natexlab{b}})]{Nagasawa87b}
M.~Nagasawa and H.~Tanaka.
\newblock On the propagation of chaos for diffusion processes with drift
  coefficients not of average form.
\newblock \emph{Tokyo J. Math.}, 10\penalty0 (2):\penalty0 403--418,
  1987{\natexlab{b}}.

\bibitem[Protter(2005)]{Protter05}
P.E. Protter.
\newblock \emph{Stochastic Integration and Differential Equations}.
\newblock Springer, Berlin, 2nd edition, 2005.

\bibitem[Rockafellar(1974)]{Rockafellar74}
R.T. Rockafellar.
\newblock \emph{Conjugate Duality and Optimization}.
\newblock SIAM, Philadelphia, 1974.

\bibitem[Sznitman(1991)]{Sznitman91}
A.-S. Sznitman.
\newblock Topics in propagation of chaos.
\newblock In P.-L. Hennequin, editor, \emph{\'Ecole d'\'Et\'e de Probabilit\'es
  de Saint Flour XIX - 1989}, volume 1464 of \emph{Lecture Notes in
  Mathematics}, pages 165--251. Springer, Berlin, 1991.

\bibitem[Talagrand(2011{\natexlab{a}})]{Talagrand11}
M.~Talagrand.
\newblock \emph{Mean Field Models for Spin Glasses. Volume I: Basic Examples}.
\newblock Springer, Berlin, 2011{\natexlab{a}}.

\bibitem[Talagrand(2011{\natexlab{b}})]{Talagrand11b}
M.~Talagrand.
\newblock \emph{Mean Field Models for Spin Glasses. Volume II: Advanced
  Replica-Symmetry and Low Temperature}.
\newblock Springer, Berlin, 2011{\natexlab{b}}.

\bibitem[Vaillancourt(1988)]{Vaillancourt88}
J.~Vaillancourt.
\newblock On the existence of random {M}c{K}ean-{V}lasov limits for triangular
  arrays of exchangeable diffusions.
\newblock \emph{Stoch. Anal. Appl.}, 6\penalty0 (4):\penalty0 431--446, 1988.

\end{thebibliography}

\end{document}